\documentclass[12pt]{article}
\usepackage{amssymb,amsmath,latexsym}
\usepackage{ mathrsfs }
\usepackage{color}
\usepackage{hyperref}
\usepackage{enumerate}
\usepackage{kotex}

\allowdisplaybreaks

\begin{document}
	\begin{doublespace}
		\def\1{{\bf 1}}
		\def\nn{\nonumber}
		\newcommand{\tG}{\tilde{\G}}
		\newcommand{\tf}{\tilde{f}}	
		\newcommand{\pinv}{\Phi^{-1}}
		\newcommand{\vinv}{\varphi^{-1}}
		\newcommand{\I}{\mathbf{1}}
		\def\sA {{\cal A}} 
		\def\sB {{\cal B}}
		\def\sC {{\cal C}}
		\def\sD {{\cal D}} \def\sE {{\cal E}} \def\sF {{\cal F}}
		\def\sG {{\cal G}} \def\sH {{\cal H}} \def\sI {{\cal I}}
		\def\sJ {{\cal J}} \def\sK {{\cal K}}
		\def\sL {{\cal L}}
		\def\sM {{\cal M}} \def\sN {{\cal N}} \def\sO {{\cal O}}
		\def\sP {{\cal P}} \def\sQ {{\cal Q}} \def\sR {{\cal R}}
		\def\sS {{\cal S}} \def\sT {{\cal T}} \def\sU {{\cal U}}
		\def\sV {{\cal V}} \def\sW {{\cal W}} \def\sX {{\cal X}}
		\def\sY {{\cal Y}} \def\sZ {{\cal Z}}
		\def\bA {{\mathbb A}} \def\bB {{\mathbb B}} \def\bC {{\mathbb C}}
		\def\bD {{\mathbb D}} \def\bE {{\mathbb E}} \def\bF {{\mathbb F}}
		\def\bG {{\mathbb G}} \def\bH {{\mathbb H}} \def\bI {{\mathbb I}}
		\def\bJ {{\mathbb J}} \def\bK {{\mathbb K}} \def\bL {{\mathbb L}}
		\def\bM {{\mathbb M}} \def\bN {{\mathbb N}} \def\bO {{\mathbb O}}
		\def\bP {{\mathbb P}} \def\bQ {{\mathbb Q}} 
		\def\R {{\mathbb R}}
		\def\bS {{\mathbb S}} \def\bT {{\mathbb T}} \def\bU {{\mathbb U}}
		\def\bV {{\mathbb V}} \def\bW {{\mathbb W}} \def\bX {{\mathbb X}}
		\def\bY {{\mathbb Y}} \def\bZ {{\mathbb Z}}
		\def\n{{\bf n}} 
		\def\Z {{\mathbb Z}}
		\def \lam{\lambda}

		\newcommand{\expr}[1]{\left( #1 \right)}
		\newcommand{\cl}[1]{\overline{#1}}
		\newtheorem{thm}{Theorem}[section]
		\newtheorem{lemma}[thm]{Lemma}
		\newtheorem{defn}[thm]{Definition}
		\newtheorem{prop}[thm]{Proposition}
		\newtheorem{corollary}[thm]{Corollary}
		\newtheorem{remark}[thm]{Remark}
		\newtheorem{example}[thm]{Example}
		\numberwithin{equation}{section}
		\def\qed{{\hfill $\Box$ \bigskip}}
		\def\NN{{\mathcal N}}
		\def\AA{{\mathcal A}}
		\def\PP{{\mathcal P}}
		\def\MM{{\mathcal M}}
		\def\BB{{\mathcal B}}
		\def\CC{{\mathcal C}}
		\def\LL{{\mathcal L}}
		\def\DD{{\mathcal D}}
		\def\FF{{\mathcal F}}
		\def\EE{{\mathcal E}}
		\def\QQ{{\mathcal Q}}
		\def\SS{{\mathcal S}}
		\def\G{{\mathscr G}}
		\def\M{{\mathscr M}}
		\def\L{{\bf L}}
		\def\K{{\bf K}}
		\def\S{{\bf S}}
		\def\A{{\bf A}}
		\def\E{{\mathbb E}}
		\def\F{{\bf F}}
		\def\P{{\mathbb P}}
		\def\N{{\mathbb N}}
		\def\eps{\varepsilon}
		\def\wh{\widehat}
		\def\wt{\widetilde}
		\def\pf{\noindent{\bf Proof.} }
		\def\pff{\noindent{\bf Proof} }
		\def\cp{\mathrm{Cap}}
		\def\kk{\mathfrak{K}}
		\title{Heat kernels of non-symmetric jump processes with exponentially decaying jumping kernel}
		
		\def\no{\noindent}

		\author{{\bf Panki Kim}
			\quad and \quad {\bf Jaehun Lee} 
		}
		\author{{\bf Panki Kim}\thanks{The research of Panki Kim  is  supported by the National Research Foundation of
Korea(NRF) grant funded by the Korea government(MSIP) (No. NRF-2015R1A4A1041675)
}
\quad and
\quad {\bf Jaehun Lee} 
\thanks{The research of Jaehun Lee is supported by
 the National Research Foundation of Korea(NRF) grant funded by the Korea governmentt(MSIP) : NRF-2016K2A9A2A13003815}
}
		
		\date{}
		
		\maketitle
		
		\begin{abstract}
In this paper we study 
the transition densities for a large class of non-symmetric Markov processes whose jumping kernels decay exponentially or subexponentially. We obtain their  upper bounds which also decay at the same rate as their jumping kernels. 
When the lower bounds of jumping kernels satisfy the weak upper scaling condition at zero, we also establish  lower bounds for the transition densities, which are sharp. 
\end{abstract}



		\section{Introduction}
		Let $d \in \N$, $\R^d$ be the $d$-dimensional Euclidian space and $\R_+=\{x\in  \R^1: x>0\}$. Define
		\begin{equation}\label{d:sL}
		\sL^\kappa f(x) := \lim_{\eps \downarrow 0} \sL^{\kappa,\eps} f(x) := \lim_{\eps \downarrow 0} \int_{\{z \in \R^d : |z| > \eps\}} \left( f(x+z) -f(x) \right) \kappa(x,z) J(|z|)dz
		\end{equation}
		where $\kappa : \R^d \times \R^d \rightarrow \R_+$ is a Borel function
		satisfying the following conditions:
		there exist positive constants $\kappa_0, \kappa_1, \kappa_2$ and $\delta \in (0,1)$ such that
		\begin{equation}\label{c:kap1}
		\kappa_0 \le \kappa(x,z) \le \kappa_1, \quad \kappa(x,z)=\kappa(x,-z) \quad \text{ for all }  x,z  \in \R^d
		\end{equation}
		and
		\begin{equation}\label{c:kap2}
		|\kappa(x,z)-\kappa(y,z)| \le \kappa_2 |x-y|^\delta \quad  \quad \text{ for all } x,y,z \in \R^d.
		\end{equation}
	The operator $\sL^\kappa$ can be regarded as the non-local
counterpart of elliptic operators in non-divergence form. In this context, the H\"older continuity of $\kappa(\cdot, z)$ in \eqref{c:kap2} is a natural assumption. 
		
		In \cite{CZ16}, Zhen-Qing Chen and Xicheng Zhang studied 
$\sL^\kappa$ and its heat kernel when  $J(r)=r^{-d-\alpha}$, $r>0$ and $\alpha \in (0,2)$. 
They proved the existence and uniqueness of the heat kernel and 
 its sharp two-sided estimates, cf.~\cite[Theorem 1.1]{CZ16} for details. 	
The methods in \cite{CZ16} are  quite robust and have been applied to non-symmetric and non-convolution operators (see \cite{BSV,CHXZ, CZ17a, CZ17c, KSV, KS17, Pe}
and references therein).
	In particular, the first named author, jointly with Renming Song and Zoran  Vondra\v{c}ek in  \cite{KSV},  studied  the operator $\LL^\kappa$ and its heat kernel
		when $J$ is comparable to jumping kernels of subordinate Brownian motions and  its L\'evy exponent satisfying a weak lower scaling condition at infinity. In this paper we consider the case that 
		 $J(r)$ decays exponentially or subexponentially when $r \to \infty$ and we obtain sharp two-sided estimates for the heat kernel of 
$\sL^\kappa$.

		Throughout this paper,  we assume $d \in \N$, and  that $J : \R_+ \rightarrow \R_+$ is continuous and non-increasing function satisfying that there exist a continuous and strictly increasing function $\phi: [0,1] \rightarrow \R_+$ with $\phi(0)=0$, and constants $b>0$, $0< \beta \le 1$ and $a\ge 1$ such that
		\begin{equation}\label{e:J1}
\frac{a^{-1}}{r^d\phi(r)} \le J(r) \le \frac{a}{r^d\phi(r)}, \quad 0<r \le 1 \quad \mbox{and} \quad J(r) \le a\exp(-br^\beta), \quad r>1.
\end{equation}
In addition, we assume that $J$ is differentiable in $\R_+$ and 
		\begin{equation}\label{e:J2}
		r \longmapsto -\frac{J'(r)}{r} \quad \mbox{ is non-increasing in } \R_+.
		\end{equation}
			Our main assumption on $\phi$ is the following weak lower scaling condition at zero: there exist $\alpha_1 \in (0,2]$ and $a_1 > 0$ such that
		\begin{equation}\label{e:lsc-phi}
		a_1 \big(\frac{R}{r}\big)^{\alpha_1} \le \frac{\phi(R)}{\phi(r)},  \quad 0< r \le R \le 1.
		\end{equation} Since we allow $\alpha_1$ to be 2, to guarantee that $J$ is to be a L\'evy density,   
		we also need the following integrability  condition for $\phi$ near zero: 
		\begin{equation}\label{e:phi}
		\int_0^1 \frac{s}{\phi(s)} ds := C_0 < \infty.
		\end{equation}
		The monotonicity of   $J(r)$  and \eqref{e:phi}  ensure the existence of an isotropic unimodal L\'evy process in $\R^d$ with the  L\'evy measure $J(|x|)dx$, which is infinite because of \eqref{e:lsc-phi} and the lower bound in \eqref{e:J1}.

Our goal is to obtain estimates of the heat kernel for $\sL^\kappa$. First we introduce the function $\G(t,x)$ which plays an important role for the estimates of heat kernel. Let us define the function $\Phi$ and $\theta$ as
\begin{equation}\label{e:dPhi}
	 \Phi(r):= \begin{cases}\displaystyle \frac{r^2}{ 2\int_0^r \frac{s}{\phi(s)}ds}, &\quad 0<r \le 1 \\
	 \Phi(1)r^2, &\quad r > 1 \end{cases}
			\end{equation} 
	and
$$
\theta(r) := \begin{cases} \displaystyle\frac{1}{r^d \Phi(r)}, &\quad r \le 1, \\
		\exp(-br^\beta)\1_{\{0<\beta<1\}} + {r^{-d-1}} \exp(-\frac{b}{5}r)\1_{\{\beta=1\}}, &\quad r > 1. \end{cases}
$$
By \eqref{e:phi}, $\int_0^r \frac{s}{\phi(s)}ds$ is integrable so that $\Phi$ is well-defined. Note that $\Phi(1)= \left( 2\int_0^1 \frac{s}{\phi(s)}ds\right)^{-1}=(2C_0)^{-1}$ is determined by $C_0$. Also, by Lemma \ref{l:wsc-Phi} we will see that $\Phi$ is a strictly increasing function in $\R_+$ and $\displaystyle \lim_{r \downarrow 0}\Phi(r)=0$, which imply that there exists an inverse function $\pinv:\R_+ \rightarrow \R_+$.  
For $t>0$ and $r>0$ define $\G(t,r)$ by
$$\G(t,r)=\G^{(d)}(t,r):= \frac{1}{t\pinv(t)^d} \land \theta(r) $$ where $a \land b := \min\{a,b\}$.
By an abuse of notation we also define 
		\begin{equation}\label{d:rho}
		\G(t,x)=\G^{(d)}(t,x) := \frac{1}{t\Phi^{-1}(t)^{d}} \land \theta(|x|), \quad t>0,\, x \in \R^d,
		\end{equation}
so $\G(t,x)=\G(t,|x|)$. Note that the definition of $\theta(r)$ for $\beta=1$ is simply technical  and it is harmless for readers to regard $\theta(r)$ as $\frac{1}{r^d\Phi(r)} \1_{\{r\le1 \}} + \exp(-\frac{b}{6}r)\1_{\{r>1  \}}$ as the upper bound of heat kernel for $\beta=1$ in  Theorems \ref{t:main1}-\ref{t:main2} below. 
		
	Let us compare $\G$ with the following function defined by
		\begin{equation}\label{d:trho}
	\tG(t,x)= \tG(t,|x|):= \frac{1}{t\pinv(t)^d} \land \frac{1}{|x|^d \Phi(|x|)}. 
	\end{equation}
By \cite[Proposition 2.1]{KSV} and  our Lemma \ref{l:usc-psi} below	we see that $\tG$ is the function used for the upper heat kernel estimate in \cite{KSV} (see Remark \ref{r:KSV} for details). It is easy to see that $\G(t,x) \le c \,  \tG(t,x)$ (see Lemma \ref{l:rho} below). Here is our main result.
		\begin{thm}\label{t:main1}
			Let $\sL^\kappa$ be the operator in \eqref{d:sL}. Assume that jumping kernel $J$ satisfies \eqref{e:J1} and \eqref{e:J2}, that $\phi$ satisfies \eqref{e:lsc-phi} and \eqref{e:phi}, and that $\kappa$ satisfies \eqref{c:kap1} and \eqref{c:kap2}. Then, there exists a unique jointly continuous function $p^\kappa(t,x,y)$ on $\R_+ \times \R^d \times \R^d$ solving
			\begin{equation}\label{e:main1-1}
			\partial_t p^\kappa(t,x,y) = \sL^\kappa p^\kappa(t,\cdot,y)(x), \quad x \neq y,
			\end{equation}
			and satisfying the following properties: 
			
			(i) (Upper bound) For every $T \ge 1$, there is a constant $c_1 >0$ such that for all $t \in (0,T]$ and $x,y \in \R^d$,
			\begin{equation}\label{e:main1-2}
			p^\kappa(t,x,y) \le c_1 t \G(t,x-y).
			\end{equation}
					(ii) (Fractional derivative) For any $x,y \in \R^d$ with $ x \neq y$, the map $t \mapsto \sL^\kappa p^\kappa(t,\cdot,y)(x)$ is continuous, and for each $T \ge 1$, there exists a constant $c_2>0$ 	such that for all $t \in (0,T]$, $\eps \in [0,1]$ and $x,y \in \R^d$,
			\begin{equation}\label{e:main1-3}
			 | \sL^{\kappa,\eps} p^\kappa(t,\cdot,y)(x) | \le c_2 \tG(t,x-y). 
			\end{equation}
			\noindent
			(iii) (\textit{Continuity}) For any bounded and uniformly continuous function $f : \R^d \rightarrow \R$,
			\begin{equation}
			\label{e:main1-4}
			\lim_{t \downarrow 0} \sup_{x \in \R^d} \left| \int_{\R^d} p^\kappa(t,x,y)f(y)dy - f(x) \right| = 0.
			\end{equation}
Furthermore, such unique function $p^\kappa(t,x,y)$ satisfies the following lower bound: for every $T \ge 1$, there exists a constant $c_3,c_4>0$ such that for all $t \in (0,T]$,
	\begin{equation}
	\label{e:main4-1}
	p^\kappa(t,x,y) \ge c_3\begin{cases} \pinv(t)^{-d}, &\quad |x-y| \le c_4 \pinv(t) \\
	t J(|x-y|), &\quad |x-y| > c_4 \pinv(t)
	\end{cases}
	\end{equation}

					The constants $c_i$, $i=1,\dots, 4$,  depend only on $d,T,a,a_1,\alpha_1, b,\beta,C_0,\delta,\kappa_0,\kappa_1$ and $\kappa_2$. 
					\end{thm}
The upper bound of the 	fractional derivative of $p^\kappa$ in  \eqref{e:main1-3}, which is a counterpart of \cite[(1.12)]{KSV}, will be  used to prove the uniqueness of heat kernel.		

		We emphasize here that unlike 	\cite[(1.21)]{KSV} we obtain \eqref{e:main4-1} without any upper weak scaling condition on $\phi$.
	The estimates in \eqref{e:main1-2}  and \eqref{e:main4-1} in Theorem \ref{t:main1} are not sharp in general. 
		However, when the jumping kernel $J$ satisfies
\begin{equation}\label{e:elowerJ} J(r) \ge a_1 \exp(-b_1 r^{\beta_1}), \quad r>1,  \end{equation}   and $\phi$ satisfies upper weak scaling condition at zero, that is, 
\begin{equation}\label{e:usc-phi}
 \frac{\phi(R)}{\phi(r)} \le a_2 \big( \frac{R}{r} \big)^{\alpha_2}, \quad 0<t \le R \le 1 
 \end{equation}
where $a_2>0$ and  $\alpha_2 \in (0,2)$, 
then the lower bound in \eqref{e:main4-1} is comparable to that in \cite[Theorem 1.2]{GKK}, which is lower heat kernel estimates for symmetric Hunt process with exponentially decaying jumping kernel. Note that $\phi$ is comparable to $\Phi$ under \eqref{e:lsc-phi} and \eqref{e:usc-phi}. Therefore, under additional assumptions \eqref{e:usc-phi} and \eqref{e:elowerJ} we have the following corollary. 

\begin{corollary}\label{c:tshke}
Let $\sL^\kappa$ be the operator in \eqref{d:sL}. Assume that jumping kernel $J$ satisfies \eqref{e:J1}, \eqref{e:J2} and \eqref{e:elowerJ}, that $\phi$ satisfies \eqref{e:lsc-phi} and \eqref{e:usc-phi}, and that $\kappa$ satisfies \eqref{c:kap1} and \eqref{c:kap2}. Then, the heat kernel $p^\kappa(t,x,y)$ for $\sL^\kappa$  satisfies the following estimates: 
for every $T \ge 1$, there is a constant $c_1 >0$ such that for all $t \in (0,T]$ and $x,y \in \R^d$,
\begin{align*}
&c_1^{-1} \left( \phi^{-1}(t)^{-d} \land  \frac{t}{|x-y|^d \phi(|x-y|)}\right)	\\
&\le p^\kappa(t,x,y) \le c_1  \left( \phi^{-1}(t)^{-d} \land  \frac{t}{|x-y|^d \phi(|x-y|)}\right), \quad |x-y| \le 1,
\end{align*}
		and
		$$
	c_1^{-1}t \exp(-b_1|x-y|^{\beta_1})	\le p^\kappa(t,x,y) \le c_1  t \theta (|x-y|), \quad |x-y| > 1.
	$$		
	The constant $c_1$ depend on $d,T,a,a_1,a_2,\alpha_1,\alpha_2,b,b_1,\beta,\beta_1,C_0,\delta,\kappa_0,\kappa_1$ and $\kappa_2$.
\end{corollary}	

Comparing to  \cite{KSV}, Corollary \ref{c:tshke} provides further precise heat kernel estimates  for the operator \eqref{d:sL} with exponential decaying function $J$. 		We remark here that, when $\beta>1$, the estimates of $p^\kappa(t,x,y)$ are different and so  the result in Corollary	\ref{c:tshke} does not hold even for symmetric L\'evy processes. See \cite{CKK1, Sz}. We will address this interesting case somewhere else. 
					
		 More properties of the heat kernel $p^\kappa(t,x,y)$ are listed in the following theorems.

	\begin{thm}\label{t:main2}
	Suppose that the assumptions of Theorem \ref{t:main1} are satisfied.

	\noindent
	(1) (Conservativeness) For all $(t,x)\in \R_+ \times \R^d$, 
	\begin{equation}\label{e:main2-1}
	\int_{\R^d}p^{\kappa}(t,x,y)\, dy =1\, .
	\end{equation}
	
	\noindent
	(2) (Chapman-Kolmogorov equation) For all $s,t > 0$ and $x,y\in \R^d$,
	\begin{equation}\label{e:main2-2}
	\int_{\R^d}p^{\kappa}(t,x,z)p^{\kappa}(s,z,y)\, dz =p^{\kappa}(t+s,x,y)\, .
	\end{equation}

	\noindent
	(3) (H\"older continuity) For every $T\ge 1$ and $\gamma \in (0,\alpha_1) \cap (0,1]$, there is a constant
	$c_1>0$ such that for all $0< t\le T$ and $x,x',y\in \R^d$ with either $x\neq y$ or $x'\neq y$,
	\begin{equation}
	\label{e:main2-3}
	|p^{\kappa}(t,x,y)-p^{\kappa}(t,x',y)|
	\le c_1 |x-x'|^{\gamma} t\, \pinv(t)^{\gamma}  (\G(t, x-y)\vee\G(t,  x'-y))\, .
	\end{equation}
The constant $c_1$ depends only on $d,T,a,a_1,\alpha_1, b,\beta,C_0,\gamma,\delta,\kappa_0,\kappa_1$ and $\kappa_2$.
\end{thm}
For $t>0$, define the operator $P_t^{\kappa}$ by
\begin{equation}\label{e:intro-semigroup}
P_t^{\kappa}f(x)=\int_{\R^d} p^{\kappa}(t,x,y)f(y)\, dy\, ,\quad x\in \R^d\, ,
\end{equation}
where $f$ is a nonnegative (or bounded) Borel function on $\R^d$, and let $P_0^{\kappa}=\mathrm{Id}$. Then
by Theorems \ref{t:main2},  $(P_t^{\kappa})_{t\ge 0}$ is a Feller semigroup with the strong Feller property.  
Let $C_b^{2, \eps} (\R^d)$ be the space of bounded twice differentiable functions in $\R^d$ whose second derivatives are uniformly H\"older continuous.	
\begin{thm}\label{t:main3}
	(1) (Generator) Let $\eps >0$. For any $f\in C_b^{2, \eps} (\R^d)$,  we have
\begin{equation}\label{e:main3-1}
\lim_{t\downarrow 0}\frac{1}{t}\left(P_t^{\kappa}f(x)-f(x)\right)=\LL^{\kappa}f(x)\, ,
\end{equation}
and the convergence is uniform.
\noindent
(2) (Analyticity) The 
semigroup $(P_t^{\kappa})_{t\ge 0}$ of $\LL^{\kappa}$ is analytic in $L^p(\R^d)$ for every $p\in [1,\infty)$.
\end{thm}

In this paper, we defined  the function $\G(t,x)$ 
from the conditions on $J$ directly, while in \cite{KSV} the function $\rho(t,x)$ is defined by the characteristic exponent of an isotropic unimodal L\'evy process with jumping kernel $J(x)dx$. The reason is that, in our situation, it is more convenient than using characteristic exponent to describe exponential decaying jumping kernel. See Remark \ref{r:KSV} below for  the connections between two definitions.

As \cite{KSV}, the approach in this paper is based on the method originally developed  in \cite{CZ16}.
In Section 2, we introduce basic setup and scaling inequalities. In addition, we obtain some convolution inequalities at Proposition \ref{l:con} in Section 2.2. The results in Section 2.2 is similar to \cite[Lemma 2.6]{KSV}, although our function $\G(t,x)$ is smaller than that in \cite{KSV}.

 In Section 3, we discuss gradient estimates for the heat kernel of isotropic unimodal L\'evy process with jumping kernel $J(|x|)dx$, which follows from the results in \cite{KS15, KS17}. We only use Proposition \ref{p:grad1} in the proof of our main theorem, but Proposition \ref{p:grad} itself is of independent interest. 
 
 In Section 4, we obtain some useful estimates on functions involving the heat kernel for the isotropic L\'evy process whose jumping kernel is $J(|x|)dx$. In Section 4.1, we improve inequalities in Propositions \ref{p:grad1} and \ref{t:p}--\ref{p:int-delta} for the symmetric L\'evy processes whose jumping kernel is $\kk(x)J(|x|)dx$, where $\kk(x)$ is symmetric and bounded between two positive constants.  As \cite[Section 3]{KSV}, we also observe the continuous dependency of the heat kernel $p^\kk$ with respect to the jumping kernel $\kk(x)J(|x|)dx$.  
 
 In Section 5, we follow the Levi's construction in \cite[Section 4]{KSV}. Note that as in \cite[Section 4]{KSV}, many results in Section 5 are derived from the estimates in Sections 2 and 4 so that we can follow \cite{KSV} for the most of proofs.
 Finally we provide the proofs of Theorems \ref{t:main1},  \ref{t:main2} and  \ref{t:main3} in Section 6.

In this paper, we use the following notations.
We will use ``$:=$" to denote a definition, which is  read as ``is defined to be".
For any two positive functions $f$ and $g$,
$f\asymp g$ means that there is a positive constant $c\geq 1$
such that $c^{-1}\, g \leq f \leq c\, g$ on their common domain of
definition. Denote {\rm diam}$(A) = \sup \{ |x-y| : x,y \in A \}$ and $\sigma(dz)= \sigma_d(dz)$ be a uniform measure in the sphere $\{ z\in \R^d : |z|=1\}$.
For a function $f:\R_+\times \R^d\to \R$, we define $f(t, x\pm z)=f(t, x+z)+f(t, x-z)$ and
\begin{equation}\label{d:del}
\delta_f(t,x;z):=f(t,x+z)+f(t,x-z)-2f(t,x) = f(t,x \pm z) - 2f(t,x).
\end{equation}

Throughout the rest of this paper, the positive constants
 $T, a,a_1, \alpha_1,b, \beta, \delta,\kappa_0,\kappa_1, \kappa_2$ and $C_i$,
$i=0,1,2,\dots $, can be regarded as fixed.
In the statements of results and the proofs, the constants $c_i=c_i(a,b,c,\ldots)$, $i=0,1,2,  \dots$, denote generic constants depending on $a, b, c, \ldots$, whose exact values are unimportant.
They start
anew in each statement and each proof.

\section{Preliminaries}
In this section we first study some elementary properties of $\Phi$ defined in \eqref{e:dPhi}.
\begin{lemma}\label{l:wsc-Phi}
	Assume that $\phi$ satisfies \eqref{e:lsc-phi} and \eqref{e:phi}. Then, $\Phi$ is continuous and strictly increasing in $(0,1]$, and satisfies
	\begin{equation}\label{e:Phi}
	\Phi(r) \le \phi(r), \quad \quad 0<r \le 1
	\end{equation}
	and
	\begin{equation}\label{e:wsc-Phi}
		a_1 \left(\frac{R}{r}\right)^{\alpha_1} \le \frac{\Phi(R)}{\Phi(r)} \le \left(\frac{R}{r}\right)^2,  \quad 0< r \le R,
	\end{equation}
	where $a_1>0$ and $\alpha_1 \in (0,2]$ are constants in \eqref{e:lsc-phi}. In particular, \eqref{e:Phi} implies $\displaystyle \lim_{r \downarrow 0} \Phi(r)=0$.
\end{lemma}
\pf Since $\phi$ is continuous in $(0,1]$, $\Phi$ is continuous in $\R_+$ by definition. Also, since $\phi$ is strictly increasing, we have
$$ \Phi(r) = \frac{r^2}{2\int_0^r \frac{s}{\phi(s)}ds} \le \frac{r^2}{2\int_0^r \frac{s}{\phi(r)}ds} = \phi(r). $$
To show that $\Phi$ is  strictly increasing,  it suffices to observe that for $0<r<1$,
$$ \left( \frac{1}{2\Phi(r)} \right)' =  \left(r^{-2} {\int_0^r \frac{s}{\phi(s)}ds } \right)' = 2r^{-3} \int_0^r s \left( \frac{1}{\phi(r)} - \frac{1}{\phi(s)} \right) ds <0.$$
Now we prove \eqref{e:wsc-Phi}. 
Clearly, by the definition of $\Phi$, \eqref{e:wsc-Phi} holds for $1 \le r \le R$. 

For $0<r \le R \le 1$, we have
${R^2}/{\Phi(R)} = \int_0^R (s/{\phi(s)})ds \ge \int_0^r ({s}/{\phi(s)})ds = {r^2}/{\Phi(r)}, $
which implies the second inequality in \eqref{e:wsc-Phi}.
Also, by change of variables and \eqref{e:lsc-phi}
\begin{align*}
\frac{\Phi(R)}{\Phi(r)}&= \frac{2\Phi(R)}{r^2} \int_0^r \frac{s}{\phi(s)}ds = \frac{2\Phi(R)}{r^2} \int_0^R \frac{(r/R) t}{ \phi((r/R)t)} (r/R)dt \\
&= \frac{2\Phi(R)}{R^2} \int_0^R \frac{t}{\phi((r/R)t)}dt \ge a_1\big(\frac{R}{r}\big)^{\alpha_1} \frac{2\Phi(R)}{R^2}\int_0^R \frac{t}{\phi(t)}dt  = a_1\big(\frac{R}{r}\big)^{\alpha_1}.
\end{align*}
For $R \ge 1 \ge r >0$, using $\Phi(R)=\Phi(1)R^2$ and above estimates we have
$$ a_1 \left( \frac{R}{r} \right)^{\alpha_1} \le a_1 \frac{ R^2}{r^{\alpha_1}} \le \frac{\Phi(R)}{\Phi(r)}= \frac{\Phi(R)}{\Phi(1)} \frac{\Phi(1)}{\Phi(r)} \le  \frac{R^2}{r^2}.$$
 \qed

Note that our main results hold for all $t \le T$, while the definition of $\G$ in \eqref{d:rho} is independent of $T$. To make our proofs simpler,  
we introduce a family of auxiliary functions which will be used mostly in proofs. 

Let $T \ge \Phi(1)$ and define $\G_T: (0,T] \times (0,\infty) \rightarrow (0,\infty)$ by
\begin{equation}\label{e:rho}
\G_T(t,r) = \begin{cases} \displaystyle \frac{1}{t\pinv(t)^d}, &\quad r \le \pinv(t), \\
\displaystyle \frac{1}{r^d\Phi(r)},  &\pinv(t)< r  \le \pinv(T), \\
 \displaystyle C_T\exp(-br^\beta)\1_{0<\beta<1} + \frac{C_T}{r^{d+1}} \exp(-\frac{b}{5}r)\1_{\beta=1}, &\quad r > \pinv(T)
\end{cases}
\end{equation} 
where $C_T:= T^{-1}\pinv(T)^{-d}\exp(b\pinv(T)^\beta)\1_{0<\beta<1 } + T^{-1}\pinv(T) \exp(\frac{b}{4}\pinv(T)) \1_{\beta=1}.$
Note that  $r \mapsto \G_T(t,r)$ is continuous and non-increasing (due to such choice of $C_T$).

Recall that $\tG(t,r)$ is defined in \eqref{d:trho}.
In the following lemma we show that $\G_T$ and  $\G(t,x)$ are comparable and less than $\tG(t,r)$. 
\begin{lemma}\label{l:rho}
(a) Let $T \ge \Phi(1)$. Then, there exists a constant $c_1=c_1(T)>0$ such that
\begin{equation}\label{e:rho_T} c_1^{-1} \G_T(t,r) \le \G(t,r) \le c_1\G_T(t,r) \end{equation}
for any $t \in (0,T]$ and $r>0$. \\
(b) There exists a constant $c_2>0$ such that 
\begin{equation}\label{e:trho} \G(t,r) \le c_2\tG(t,r).  \end{equation}
for any $t>0$ and $r>0$. The constant $c_1$ depends on $d,b,T,\pinv(T),\beta$ and $C_0$, and $c_2$ depends on $d,b,\beta$ and $C_0$.
\end{lemma}

\pf (a) Define 
\begin{equation*}
\theta_T(r) := \begin{cases} r^{-d}\Phi(r)^{-1}, &\quad r \le \pinv(T), \\
C_T\exp(-br^\beta)\1_{0<\beta<1 } + {C_T}{r^{-d-1}}\exp(-\frac{b}{5}r)\1_{\beta=1}, &\quad r > \pinv(T). \end{cases}
\end{equation*}
 Note that $r \mapsto \theta_T(r)$ is strictly decreasing and  $\theta_T(\pinv(t))=\frac{1}{t\pinv(t)^d}$ for any $0<t \le T$. Thus we can obtain
\begin{equation}\label{e:rel}
\theta_T(r) \le \frac{1}{t\pinv(t)^d}\quad \mbox{if and only if} \quad t \le \Phi(r).
\end{equation}
By \eqref{e:rho} and \eqref{e:rel} we have
\begin{equation}\label{e:rel1}  \G_T(t,r) =\frac{1}{t \pinv(t)^d}  \land \theta_T(r). \end{equation}
Let
$$\displaystyle  M_T:= \begin{cases}\sup_{1 \le r \le \pinv(T)} \frac{1}{r^d\Phi(r)} \exp(b r^\beta)  &\quad \mbox{for} \quad 0<\beta<1, \\ \sup_{1 \le r \le \pinv(T)} \frac{r}{\Phi(r)} \exp(\frac{b}{5} r) &\quad \mbox{for} \quad \beta=1 \end{cases}$$
and
$$\displaystyle m_T:=  \begin{cases}\inf_{1 \le r \le \pinv(T)}\frac{1}{r^d\Phi(r)} \exp(b r^\beta) &\quad \mbox{for} \quad 0<\beta<1, \\ \inf_{1 \le r \le \pinv(T)} \frac{r}{\Phi(r)} \exp(\frac{b}{5} r) &\quad \mbox{for} \quad \beta=1 \end{cases}$$
Then, for $0<\beta <1$,
$$  \theta(r)=\begin{cases} \frac{1}{r^d\Phi(r)} = \theta_T(r), &\quad r \le 1, \\
\exp(-br^\beta) \ge M_T^{-1} \frac{1}{r^d\Phi(r)} = M_T^{-1}\theta_T(r), &\quad 1< r \le \pinv(T), \\
\exp(-br^\beta) \le m_T^{-1}\theta_T(r), &\quad 1< r \le \pinv(T), \\
\exp(-br^\beta) = C_T^{-1} \theta_T(r), & \quad r>\pinv(T)
\end{cases}  $$
and for $\beta=1$,
$$  \theta(r)=\begin{cases} \frac{1}{r^d\Phi(r)} = \theta_T(r), &\quad r \le 1, \\
\frac{1}{r^{d+1}}\exp(-\frac{b}{5}r) \ge M_T^{-1} \frac{1}{r^{d}\Phi(r)} =  M_T^{-1}\theta_T(r), &\quad 1< r \le \pinv(T), \\
\exp(-br^\beta) \le m_T^{-1}\theta_T(r), &\quad 1< r \le \pinv(T), \\
\frac{1}{r^{d+1}}\exp(-\frac{b}{5}r) = C_T^{-1} \theta_T(r) , & \quad r>\pinv(T).
\end{cases}  $$
Thus, for any $0<\beta \le 1$ and $r>0$,
$$
(1 \wedge M_T^{-1} \wedge C_T^{-1})\theta_T(r) \le \theta(r) \le(1 \lor m_T^{-1} \lor C_T^{-1}) \theta_T(r).$$
Using this and \eqref{e:rel1} we arrive \eqref{e:rho_T}.

\noindent
(b) Clearly we have  $\tG(t,r) = \G(t,r)$ for $r \le 1$. For any $r > 1$ and $0<\beta <1$ we have
$$ \tG(t,r)= \frac{1}{r^d\Phi(r)}  \ge \left(\sup_{s \ge 1}  s^d\Phi(s) \exp(-bs^\beta)\right)^{-1} \exp(-br^\beta)= c(\beta) \G(t,r).$$ 
Similarly, for $r>1$ and $\beta=1$
$$ \tG(t,r)= \frac{1}{r^{d+2}\Phi(1)} \ge \left(\sup_{s \ge 1}  \frac{\Phi(s)}{s} \exp(-\frac{b}{5}s)\right)^{-1} \frac{1}{r^{d+1}} \exp(-\frac{b}{5}r)= c(1)\G(t,r). $$
Combining above estimates with \eqref{e:rho_T} we arrive \eqref{e:trho} with $c_2 = c(\beta) \land c_1^{-1}$. \qed

In the following remark we will see that our $\tG(t,x)$ and the function $\rho(t,x)$ in \cite{KSV} are comparable.
\begin{remark}\label{r:KSV} \rm
	Let ${r(t,r)}:= \psi^{-1}(t^{-1})^d \land [t {\psi(r^{-1})}r^{-d}]$ as in \cite{KSV}, where $\psi$ is the characteristic exponent with respect to the L\'evy process whose jumping kernel is $J(|y|)dy$. By Lemma \ref{l:usc-psi} below we have $\psi(r^{-1})^{-1} \asymp \Phi(r)$ for all $r>0$, which implies that ${r(t,r)}/{t} \asymp \tG(t,r)$ for all $r>0$.  Thus, by \cite[Proposition 2.1]{KSV} we conclude that $\tG(t,x)$ is comparable to the function $\rho(t,x)$ in \cite{KSV}.
\end{remark}

\subsection{Basic scaling inequalities.} 
We start with weak scaling condition for the inverse ) of $\Phi$. In this subsection we assume that $\phi$ satisfies \eqref{e:lsc-phi}.
\begin{lemma}\label{l:wsc-inv}
For any $0<r \le R$,
		\begin{equation}\label{e:wsc-inv}
	 \big(\frac{R}{r}\big)^{1/2}  \le \frac{\pinv(R)}{\pinv(r)} \le a_1^{-1/\alpha_1}  \big(\frac{R}{r}\big)^{1/\alpha_1}
		\end{equation}
	where $a_1$ and $\alpha_1$ are constants in \eqref{e:lsc-phi}. 
\end{lemma}
\pf Letting $(r,R)=(\pinv(r),\pinv(R))$ in \eqref{e:wsc-Phi}, we have that for $0<r \le R$,
	$$a_1 \big( \frac{\pinv(R)}{\pinv(r)} \big)^{\alpha_1} \le  \frac{R}{r} = \frac{\Phi(\pinv(R))}{\Phi(\pinv(r))} \le \big( \frac{\pinv(R)}{\pinv(r)} \big)^{2},$$
	which implies \eqref{e:wsc-inv}.	\qed
	
Now we introduce some scaling properties of $\G$ which will be used throughout this paper.
\begin{lemma}\label{l:rho-scaling}
	Let $T\ge 1$ and $0<\eps$. Then, there exist constants $c_1,c_2>0$ such that for any $0<t\le T$, $x \in \R^d$ and $z \in \R^d$ satisfying $\Phi(|z|) \le t$,
	\begin{equation}\label{e:rho-scaling(a)}
	\G(\eps t, x) \le c_1\G(t,x)
	\end{equation}
	and
	\begin{equation}\label{e:rho-scaling(b)}
	\G(t,x+z) \le c_2\G(t,x),
	\end{equation}
	where $c_1$ depends only on $d, a_1,\alpha_1,\eps$, and $c_2$  depends only on $d,T, a_1,\alpha_1,b,\beta$ and $C_0$.		
\end{lemma}
\pf (a) Since $t \mapsto \G(t,x)$ is non-increasing, we can assume $\eps <1$ without loss of generality. By \eqref{e:wsc-inv}, there is a constant $c_1(\eps)>1$ satisfying
$$ \frac{1}{\eps t \pinv(\eps t)^d} \le \frac{c_1}{t \pinv(t)^d} $$
for $t \le T$. Thus, we arrive
\begin{align*} 
\G(\eps t,x)=\frac{1}{\eps t \pinv(\eps t)^d}  \wedge \theta(|x|)  \le \frac{c_1}{t\pinv(t)^d} \land \theta(|x|) \le c_1\G(t,x). 
\end{align*}
\noindent 
(b) We claim that \eqref{e:rho-scaling(b)} holds with the ) $\G_T(t,x)$. In other words, there exists a constant $c_2>0$ such that
$$ \G_T(t,x+z) \le c_2 \G_T(t,x), \quad x \in \R^d, \,\, \Phi(|z|) \le t. $$
Since $r \mapsto \G_T(t,r)$ is non-increasing, it suffices to show that there exists $c_3>0$ such that for any $0<t \le T$ and $r>0$,
\begin{equation}\label{e:scale}
\G_T(t,r) \le c_3 \exp(b\pinv(t)^\beta)\G_T(t,r+\pinv(t)).
\end{equation}
Indeed, since $0<t \le T$, \eqref{e:scale} implies our claim with $c_2=c_3\exp(b\pinv(T)^\beta)$.  

We prove \eqref{e:scale} by considering several cases separately. Firstly when $r \le \pinv(T) - \pinv(t)$, using \eqref{e:rho} we have
\begin{align*}
\G_T(t,r+\pinv(t)) &= \frac{1}{(r+\pinv(t))^d \Phi(r+\pinv(t))} \\ &\ge  \frac{1}{(2\pinv(t))^d \Phi(2\pinv(t))} \land \frac{1}{(2r)^d \Phi(2r)} \\ &\ge c_4 \left(\frac{1}{t\pinv(t)^d} \land \frac{1}{r^d\Phi(r)} \right) = c_4\G_T(t,r),
\end{align*}
The second line above follows from \eqref{e:wsc-Phi}.

When $r \ge \pinv(T)$ and $0<\beta<1$, using \eqref{e:rho} and triangular inequality $r^\beta + \pinv(t)^\beta \ge (r+\pinv(t))^\beta$ we get
\begin{align*}
\G_T(t,r+\pinv(t)) = \exp(-b(r+\pinv(t))^\beta) \ge \exp(-b\pinv(t)^\beta-br^\beta)= \exp(-b\pinv(t)^\beta) \G_T(t,r). 
\end{align*}
Similarly, for $r \ge \pinv(T)$ and $\beta=1$ we have
\begin{align*}
\G_T(t,r+\pinv(t))&= C_T \frac{1}{(r+\pinv(t))^{d+1}}\exp(-\frac{b}{5} (r+\pinv(t))) \ge C_T \frac{1}{(2r)^{d+1}}\exp(-\frac{b}{5}\pinv(t)- \frac{b}{5}r) \\&= 2^{-d-1} \exp(-\frac{b}{5} \pinv(t)) \G_T(t,r) \ge 2^{-d-1} \exp(-b\pinv(t))\G_T(t,r).
\end{align*}
When $\pinv(T)-\pinv(t) \le r \le \pinv(T)$, combining above estimates we arrive
\begin{equation*}
\G_T(t,r+\pinv(t)) \ge 2^{-d-1}\exp(-b\pinv(t)^\beta)\G_T(t,\pinv(T)) \ge c_5 \exp(-b\pinv(t)^\beta)\G_T(t,r).
\end{equation*}
Note that $r \mapsto \G_T(t,r)$ is continuous at $r=\pinv(T)$. Therefore, we conclude \eqref{e:scale}. Applying \eqref{e:rho_T} for \eqref{e:scale} we arrive our desired estimate \eqref{e:rho-scaling(b)}. \qed

\subsection{Convolution inequalities}

In this section, we obtain some convolution inequalities for $\G(t,x)$ which will be used for Levi's method in Section 5. To get these inequalities we will use some estimates in \cite[Section 2]{KSV}. Note that by Remark \ref{r:KSV} we already have convolution inequalities for $\tG(t,x)$ (e.g. \cite[Proposition 2.8]{KSV}). For $a,b>0$, let $B(a,b)= \int_0^1 s^{a-1} (1-s)^{b-1}ds=\frac{(a+b-1)!}{(a-1)!(b-1)!}$ be the beta function. 

Using \eqref{e:wsc-inv}, the proof of the following lemma is same as the one in \cite[Lemma 2.3]{KSV}. Thus we skip the proof. 

\begin{lemma}\label{l:con-phi}
	Assume that $\phi$ satisfies \eqref{e:lsc-phi} and $\gamma, \delta \ge 0$, $\eta, \theta\in \R$ are constants satisfying 	${\bf 1}_{\gamma\ge 0}(\gamma/2)  +{\bf 1}_{\gamma< 0}(\gamma/\alpha_1)+\delta/2+1-\eta>0$. Then for every $t >0$,  we have
	\begin{equation}\label{e:con-phi}
	\int_0^t s^{-\eta} \pinv(s)^{\gamma} (t-s)^{-\theta}\pinv(t-s)^{\delta}\, ds \le B(\delta/2+1-\theta, \gamma/2+1-\eta)t^{1-\eta-\theta}\pinv(t)^{\gamma+\delta}\,.
	\end{equation}
\end{lemma}

For $0 \le s \le t$, let $g(s) := t^\beta + (2^\beta-1)s^\beta - (t+s)^\beta$. Then we can easily check that $g(0) = g(t) =0$ and $$g'(s)= \beta \left( (2^\beta-1)s^{\beta-1} - (t+s)^{\beta-1} \right) \begin{cases}  \,\,\ge 0, &\quad s \in [0, kt ], \\
\,\,\le 0, &\quad s \in [kt,t],  \end{cases}$$
where $k:= ((2^\beta -1)^{\frac{1}{\beta-1}}-1)^{-1}\in (0,1)$ is the constant satisfying $g'(kt)=0$.
Thus, we conclude that $g(s) \ge 0$ for any $0 \le s \le t$, which implies 
\begin{equation}\label{e:exp2}
 t^\beta + s^\beta \ge (t+s)^\beta + (2-2^\beta)  (t^\beta \land s^\beta), \quad \text{ for all }0<\beta <1 \text{ and }t,s>0.
\end{equation}

Using \eqref{e:exp2} we prove the following lemma, which we need for our convolution inequalities.
\begin{lemma}\label{l:exp}
	(a) Let  $0<\beta<1$ and $b>0$. Then, there exists a constant $c_1>0$ such that for any $x \in \R^d$,
	\begin{equation}\label{e:exp}
	\int_{\R^d} \exp(-b|x-z|^\beta -b|z|^\beta) dz \le c_1 \exp(-b |x|^\beta).
	\end{equation}
	(b) There exists a constant $c_2>0$ such that for any $x \in \R^d$ with $|x| \ge 1$,
	\begin{equation}
	\label{e:exp1}
	\int_{\R^d} ( |x-z|^{-d-1} \land 1 )( |z|^{-d-1} \land 1 ) dz \le c_2 |x|^{-d-1}.
	\end{equation}
	The constant $c_1$ depends only on $b$, $d$ and $\beta$, and $c_2$ depends only on $d$.
\end{lemma}
\pf (a)  Let 
$$ c_1=2 \int_{\R^d} \exp(-b(2-2^\beta)|z|^\beta)dz <\infty. $$
Using \eqref{e:exp2} for the second line, we arrive 
\begin{align*}
\int_{\R^d} \exp(-&b|x-z|^\beta - b|z|^\beta) dz \le \int_{\R^d} \exp(-b|x|^\beta) \exp \big(-b(2-2^\beta) (|z|^\beta \land |x-z|^\beta)\big)dz \\
&\le \exp(-b|x|^\beta) \left( \int_{\R^d} \exp(-b(2-2^\beta)|z|^\beta)dz + \int_{\R^d} \exp(-b(2-2^\beta)|x-z|^\beta)dz  \right) \\ &=   c_1\exp(-b|x|^\beta).
\end{align*}
 This proves \eqref{e:exp}.

\noindent (b) Using $|x-z|^{-1} \land |z|^{-1} \le 2|x|^{-1}$, we have
\begin{align*}
&\int_{\R^d} ( |x-z|^{-d-1} \land 1 )( |z|^{-d-1} \land 1 ) dz\\
 \le&  (\frac{2}{|x|})^{d+1} \left(\int_{|x-z| \ge |z|} (|z|^{-d-1} \land 1) dz + \int_{|x-z| < |z|} (|x-z|^{-d-1} \land 1) dz\right)\\ \le& (\frac{2}{|x|})^{d+1} \left(\int_{\R^d} (|z|^{-d-1} \land 1) dz + \int_{\R^d} (|x-z|^{-d-1} \land 1) dz\right)  :=	c_2 |x|^{-d-1}.
\end{align*}
This concludes the lemma. \qed

 For $\gamma, \delta \in \R$, $t>0$ and $x \in \R^d$ we define
$$
\G_{\gamma}^{\delta}(t,x):=\pinv(t)^{\gamma}(|x|^{\delta}\wedge 1) \G(t,x) \quad \mbox{and} \quad \tG_{\gamma}^\delta(t,x) := \pinv(t)^{\gamma} (|x|^{\delta}\wedge 1) \tG(t,x).
$$
Note that 
$\G_0^0(t,x)=\G(t,x)$, and  $\tG_{\gamma}^\delta(t,x)$ is comparable to the function $\rho_{\gamma}^\delta(t,x)$ in \cite{KSV} by Remark \ref{r:KSV}.
Also, we can easily check that  for $T \ge \Phi(1)$, 
\begin{align}
&\G_{\gamma_1}^\delta(t,x) \le \pinv(T)^{\gamma_1-\gamma_2}\G_{\gamma_2}^\delta(t,x),  \qquad \,\,\,\, (t,x) \in (0,T] \times \R^d, \quad \, \gamma_2 \le \gamma_1, \label{e:rho-gam}\\ 
&\G_{\gamma}^{\delta_1}(t,x) \le \G_{\gamma}^{\delta_2}(t,x), \quad \, \, \, \, \qquad \qquad \qquad (t,x) \in (0,\infty) \times \R^d, \quad 0 \le \delta_2 \le \delta_1. \label{e:rho-del}
\end{align}
We record the following inequality which immediately follows   from \eqref{e:rho-gam} and \eqref{e:rho-del}: for any $T \ge \Phi(1)$, $\delta \ge 0$ and $(t,x) \in (0,T] \times \R^d$, 
\begin{equation}\label{e:rho1}
\big(\G_0^\delta + \G_\delta^0\big) (t,x) \le \big(\pinv(T)^{\delta} +1 \big) \G(t,x) \le 2 \pinv(T)^{\delta} \G(t,x).
\end{equation}

Now we are ready to introduce convolution inequalities for $\G(t,x)$. 

\begin{prop}\label{l:con} Assume that $\phi$ satisfies \eqref{e:lsc-phi}. Let  $T \ge 1$ and $0<\alpha<\alpha_1$. \\
	(a) There exists a constant $c=c(d,T, a_1,\alpha, \alpha_1)>0$ such that for any
	$0< t\le T$, 
	 $\delta\in [0,\alpha]$ and $\gamma \in \R$,
	\begin{equation}\label{e:con(a)}
	\int_{\R^d}\tG_{\gamma}^{\delta}(t,x)\, dx \le c t^{-1}\pinv(t)^{\gamma+\delta}\, .
	\end{equation}
	(b)  There exists $C= C(\alpha,T)=C(d,T,a_1,\alpha, \alpha_1, b,\beta)> 0$ such that for all $x \in \R^d, \delta_1,\delta_2  \ge 0$ with $\delta_1+ \delta_2  \le \alpha$,
	$\gamma_1,\gamma_2\in \R$ and  
	$0<s<t\le T$, 
	\begin{align}\label{e:con(b)}
	\int_{\R^d}\G_{\gamma_1}^{\delta_1}(t-s,x-z) \G_{\gamma_2}^{\delta_2}(s,z)\, dz 
	\le & C\Big((t-s)^{-1}\pinv(t-s)^{\gamma_1+\delta_1+\delta_2}\pinv(s)^{\gamma_2}\G(t,x)\nn\\
	&  +\pinv(t-s)^{\gamma_1}s^{-1}\pinv(s)^{\gamma_2+\delta_1+\delta_2}  \G(t,x)\nonumber \\
	& +(t-s)^{-1}\pinv(t-s)^{\gamma_1+\delta_1}\pinv(s)^{\gamma_2}\G_0^{\delta_2}(t,x)\nonumber \\
	& +\pinv(t-s)^{\gamma_1}s^{-1}\pinv(s)^{\gamma_2+\delta_2}\G_0^{\delta_1}(t,x)\,\Big).
	\end{align}
	In particular, letting $\gamma_1=\gamma_2=\delta_1=\delta_2=0$ in \eqref{e:con(b)} we have
	\begin{equation}\label{e:con}
	\int_{\R^d} \G(t-s,x-z)\G(s,z)dz \le 2C \big( s^{-1} + (t-s)^{-1} \big) \G(t,x).
	\end{equation}
	\noindent
	(c) For all $x \in \R^d$, $0<t \le T$, $\delta_1,\delta_2  \ge 0$ and $\theta,\eta \in[0,1]$ satisfying $\delta_1+ \delta_2  \le \alpha$, 	${\bf 1}_{\gamma_1\ge 0}(\gamma_1/2)  +{\bf 1}_{\gamma_1< 0}(\gamma_1/\alpha_1)
	+\delta_1/2+1-\theta>0$ and
	${\bf 1}_{\gamma_2\ge 0}(\gamma_2/2)  +{\bf 1}_{\gamma_2< 0}(\gamma_2/\alpha_1)+\delta_2/2+1-\eta>0$, we have a constant $C_2>0$ satisfying
	\begin{align}\label{e:con(c)}
	\lefteqn{\int_0^t \int_{\R^d} (t-s)^{1-\theta}\G_{\gamma_1}^{\delta_1}(t-s, x-z)s^{1-\eta} \G_{\gamma_2}^{\delta_2}(s,z)\, dz\, ds }\nonumber \\
	\le &C_2 t^{2-\theta-\eta}\left(\G_{\gamma_1+\gamma_2+\delta_1+\delta_2}^0+\G_{\gamma_1+\gamma_2+\delta_2}^{\delta_1}+\G_{\gamma_1+\gamma_2+\delta_1}^{\delta_2}\right)(t,x)\, 
	\end{align}
	for any $0<t \le T$ and $x \in \R^d$. Moreover, when $\gamma_1, \gamma_2 \ge 0$ we further have
	\begin{align}\label{e:con(d)}
	C_2=4C \,B\left(\frac{\gamma_1+\delta_1}{2}+1-\theta,\frac{\gamma_2+\delta_2}{2}+1-\eta\right).
	\end{align}
\end{prop}
\pf (a) See  \cite[Lemma 2.6(a)]{KSV}. \\ 
(b) By \eqref{e:rho_T}, it suffices to show \eqref{e:con(b)} with the function $(\G_T)^\delta_\gamma(t,x):=\pinv(t)^{\gamma}(|x|^{\delta}\wedge 1) \G_T(t,x).$ Without loss of generality we assume $T \ge \Phi(1)$ and for notational  convenience we drop $T$ in the notations 
so we  use $\G(t,x)$ and $\G^\delta_\gamma(t,x)$ instead of $\G_T(t,x)$ and $(\G_T)^\delta_\gamma(t,x)$ respectively.

First let $|x| \le \pinv(T)$. By Remark \ref{r:KSV} and \cite[Lemma 2.6(b)]{KSV}, we already have that there exists $c_1>0$ satisfying
	\begin{align*}
\int_{\R^d}\tG_{\gamma_1}^{\delta_1}(t-s,x-z) \tG_{\gamma_2}^{\delta_2}(s,z)\, dz 
\le & c_1\Big((t-s)^{-1}\pinv(t-s)^{\gamma_1+\delta_1+\delta_2}\pinv(s)^{\gamma_2}\tG(t,x)\\
&  +\pinv(t-s)^{\gamma_1}s^{-1}\pinv(s)^{\gamma_2+\delta_1+\delta_2}  \tG(t,x) \\
& +(t-s)^{-1}\pinv(t-s)^{\gamma_1+\delta_1}\pinv(s)^{\gamma_2}\tG_0^{\delta_2}(t,x) \\
& +\pinv(t-s)^{\gamma_1}s^{-1}\pinv(s)^{\gamma_2+\delta_2}\tG_0^{\delta_1}(t,x)\,\Big).
\end{align*}
Note that $\G(t,x)=\tG(t,x)$ by \eqref{e:rho} since $|x| \le \pinv(T)$. Using \eqref{e:trho} for the left-hand side and $\G(t,x) = \tG(t,x)$ for the right-hand side, we obtain \eqref{e:con(b)} for $|x| \le \pinv(T)$. \\

Now assume $|x| > \pinv(T)$ and observe that
\begin{align*}
&\int_{\R^d} \G_{\gamma_1}^{\delta_1}(t-s,x-z) \G_{\gamma_2}^{\delta_2} (s,z) dz \\ &= \left(\int_{\underset{|x-z| > \pinv(T)}{  |z| > \pinv(T),}} + \int_{\underset{|x-z| \le \pinv(T)}{  |z| > \pinv(T),}} + \int_{\underset{|x-z| > \pinv(T)}{  |z| \le \pinv(T),}} + \int_{\underset{|x-z| \le \pinv(T)}{  |z| \le \pinv(T),}} \right)\G_{\gamma_1}^{\delta_1}(t-s,x-z) \G_{\gamma_2}^{\delta_2} (s,z) dz \\
&:= I_1 + I_2 + I_3 + I_4.
\end{align*}
First we assume $0<\beta<1$ and obtain upper bounds for $I_i$, $i=1, \dots 4$.
For $I_1$, using $\pinv(T) \ge 1$ we have
\begin{align}
I_1 &=\int_{|x-z| > \pinv(T), |z| > \pinv(T)} \G_{\gamma_1}^{\delta_1}(t-s,x-z) \G_{\gamma_2}^{\delta_2} (s,z) dz \nn \\
&= \int_{|x-z| > \pinv(T), |z| > \pinv(T)} \pinv(t-s)^{\gamma_1} \big(|x-z|^{\delta_1} \land 1 \big) \G(t-s,x-z) \pinv(s)^{\gamma_2} \big( |z|^{\delta_2} \land 1 \big) \G (s,z) dz \nn \\
&= \pinv(t-s)^{\gamma_1} \pinv(s)^{\gamma_2} \int_{|x-z| > \pinv(T), |z| > \pinv(T)} \exp(-b|x-z|^\beta - b|z|^\beta)dz. \label{2}
\end{align}
By \eqref{e:exp} we obtain
\begin{align*}
I_1 &\le c_1 \pinv(t-s)^{\gamma_1} \pinv(s)^{\gamma_2} \exp(-b|x|^\beta) = c_1 \pinv(t-s)^{\gamma_1} \pinv(s)^{\gamma_2} \G(t,x) \\
&\le c_2 (t-s)^{-1}\pinv(t-s)^{\gamma_1 + \delta_1 + \delta_2} \pinv(s)^{\gamma_2} \G(t,x).
\end{align*}
where we used $\delta_1, \delta_2 \ge 0$ and $t-s \le T$ for the last line. For the estimates of $I_2,I_3$ and $I_4$ we omit  counterpart of the last line above.\\

For $I_2$, using \eqref{e:rho} we have
\begin{align*}
I_2&=\int_{|x-z| \le \pinv(T), |z| > \pinv(T)} \G_{\gamma_1}^{\delta_1}(t-s,x-z) \G_{\gamma_2}^{\delta_2} (s,z) dz \\
&=\pinv(t-s)^{\gamma_1} \pinv(s)^{\gamma_2} \int_{|x-z| \le \pinv(T), |z| > \pinv(T)} \tG^{\delta_1}_0(t-s,x-z) \exp(-b|z|^\beta) dz.
\end{align*}
Since $|x-z| \le \pinv(T)$, using triangular inequality we have \begin{equation}\label{1}
\exp(-b|z|^\beta) \le \exp(-b|x|^\beta) \exp(b|x-z|^\beta) \le \exp(b\pinv(T)^\beta) \exp(-b|x|^\beta).
\end{equation} Thus by \eqref{e:con(a)},
\begin{align*}
I_2 &\le \pinv(t-s)^{\gamma_1} \pinv(s)^{\gamma_2} \exp(-b|x|^\beta) \int_{\R^d} \tG_0^{\delta_1} (t-s,x-z)dz \\ &\le c_3 (t-s)^{-1}\pinv(t-s)^{\gamma_1+\delta_1} \pinv(s)^{\gamma_2} \G(t,x).
\end{align*}
By the similar way, we obtain
$$ I_3 \le c_3 s^{-1}\pinv(t-s)^{\gamma_1} \pinv(s)^{\gamma_2+\delta_2} \G(t,x). $$
When $|x| \ge 2\pinv(T)$, we have $I_4=0$. So we can assume $|x| < 2\pinv(T)$ without loss of generality for the estimate of $I_4$. By \eqref{e:trho} We have
\begin{align*}
I_4 \le \int_{\R^d} \tG_{\gamma_1}^{\delta_1}(t-s,x-z) \tG_{\gamma_2}^{\delta_2}(s,z) dz \le c_4 \tG(t,x).
\end{align*}
Using $\tG(t,x) \le \tG(t,\pinv(T)) = \G(t,\pinv(T)) \le \exp(b\pinv(T)^\beta) \G(t,2\pinv(T)) \le \exp(b \pinv(T)^\beta) \G(t,x)$, we can obtain desired estimates. 
Combining estimates for $I_1,I_2,I_3$ and $I_4$, we arrive \eqref{e:con(b)} for $0<\beta<1$.  \\

For the case $\beta=1$, estimate for $I_4$ is same as above. For $I_2$ and $I_3$, instead of  \eqref{1} we argue as the following: using $|x-z| \le \pinv(T)$ and $|x|,|z| \ge \pinv(T)$, we have
$$\frac{1}{|z|^{d+1}}\exp(-\frac{b}{5}|z|) \le \frac{2^{d+1}}{|x|^{d+1}}\exp(\frac{b}{5}\pinv(T))\exp(-\frac{b}{5}|x|).$$
For $I_1$, following \eqref{2} and using \eqref{e:exp1} for the fourth line and \eqref{e:wsc-inv} for the fifth  line we have
\begin{align*}
I_1 &= \pinv(t-s)^{\gamma_1} \pinv(s)^{\gamma_2} \int_{|x-z| > \pinv(T), |z| > \pinv(T)} \frac{1}{|x-z|^{d+1} |z|^{d+1}}\exp(-\frac{b}{5}|x-z| - \frac{b}{5}|z|)dz \\
&\le c_1 \pinv(t-s)^{\gamma_1} \pinv(s)^{\gamma_2} \exp(-\frac{b}{5}|x|) \int_{|x-z| > 1, |z| > 1} \frac{1}{|x-z|^{d+1} |z|^{d+1}} dz
\\ &\le c_1 \pinv(t-s)^{\gamma_1} \pinv(s)^{\gamma_2} \exp(-\frac{b}{5}|x|) \int_{\R^d} \big( 1 \land |x-z|^{-d-1} \big) \big( 1 \land |z|^{-d-1} \big) dz \\
&\le  c_2 \pinv(t-s)^{\gamma_1} \pinv(s)^{\gamma_2} \frac{1}{|x|^{d+1}}\exp(-\frac{b}{5}|x|) = c_2 \pinv(t-s)^{\gamma_1} \pinv(s)^{\gamma_2} \G(t,x) \\
&\le c_3 (t-s)^{-1}\pinv(t-s)^{\gamma_1 + \delta_1 + \delta_2} \pinv(s)^{\gamma_2} \G(t,x).
\end{align*}

\noindent (c) Integrating \eqref{e:con(b)} with respect to $s$ from $0$ to $t$. With \eqref{e:con-phi}, we can follow the proof of \cite[Lemma 2.6(c)]{KSV}. \qed

\section{Heat kernel estimates for L\'evy processes} 
Following the framework of \cite{CZ16, KSV}, we need estimates of derivatives of the heat kernel for the symmetric L\'evy process whose jumping kernel is $J(|y|)$ (see, for example, \cite[Proposition 3.2]{KSV}). To be more precise, in our case, 
to get the upper bound of heat kernel for non-symmetric operator of the form \eqref{d:sL}, 
we need correct upper bounds  of the first and second order derivatives of 
the heat kernel 
  for unimodal L\'evy processes.  In this section, we will prove that \eqref{e:J1} and \eqref{e:J2} are sufficient condition for the estimates of the second order derivatives in Proposition \ref{p:grad1}, which decay exponentially {or subexponentially}.

\subsection{Settings}
In this section, we fix $T \le [1, \infty)$ and let $\nu(dy)=\nu(|y|)dy$ be an isotropic measure in $\R^d$ satisfying $\int_{\R^d} \big(1 \land |y|^2 \big)\nu(dy)<\infty$. Throughout this section we further assume that $\nu : \R_+ \rightarrow \R_+$ is non-increasing, differentiable function.

Here are our goals in this section.

\begin{prop}\label{p:grad}
	Let $X$ be an isotropic unimodal L\'evy process in $\R^d$ with L\'evy measure $\nu(|y|)dy$ satisfying the following assumptions: $\phi$ is a nondecreasing function with $\phi(0)=0$ satisfying \eqref{e:lsc-phi} and \eqref{e:phi}, and there exist constants $a>0$ and $0<\beta \le 1$ such that
	\begin{equation}\label{c:nu}
	\frac{a^{-1}}{r^d\phi(r)} \le \nu(r) \le \frac{a}{r^d\phi(r)}, \quad 0<r \le 1 \quad \mbox{and} \quad \nu(r) \le a \exp(-br^\beta), \quad r>1.
	\end{equation}
 Then its transition density $x \mapsto p_t(x)$ is in $C^\infty_b(\R^d)$  and satisfies gradient estimates
	\begin{equation}\label{e:grad}
	|\nabla^k_x p_t(x)| \le c t\,  
	\G_{-k}^0(t,x)=
	\pinv(t)^{-k} \left(\frac{1}{t\Phi^{-1}(t)^{d}} \land \theta(|x|)	\right)
	, \quad k=0,1
	\end{equation}
	for any $0<t \le T$ and $x \in \R^d$. The constant $c$ depends only on $k,d,T,a,a_1,\alpha_1,b,\beta$ and $C_0$.
\end{prop}

With the above result, we can obtain the second gradient estimate for the isotropic unimodal L\'evy process whose jumping kernel satisfies \eqref{e:J1} and \eqref{e:J2}.
\begin{prop}\label{p:grad1}
	Suppose that $\phi$ is a nondecreasing function with $\phi(0)=0$ satisfying \eqref{e:lsc-phi} and \eqref{e:phi}, and  that L\'evy measure $J(|y|)dy$ satisfies \eqref{e:J1} and \eqref{e:J2} with $0<\beta \le 1$. Then, its corresponding transition density $x \mapsto p(t,x)$ is in $C^\infty_b(\R^d)$ and satisfies gradient estimates
	\begin{equation}\label{e:grad2}
	|\nabla^k_x p(t,x)| \le c t\,\G_{-k}^0(t,x)=
	\pinv(t)^{-k} \left(\frac{1}{t\Phi^{-1}(t)^{d}} \land \theta(|x|)	\right), \quad k=0,1,2
	\end{equation}
	for any $0<t \le T$ and $x \in \R^d$. The constant $c $ depends only on $k,d,T,a,a_1,\alpha_1,b,\beta$ and $C_0$.
\end{prop}

In the next subsection, we prove Propositions \ref{p:grad} and \ref{p:grad1}.
\subsection{Proof of Propositions \ref{p:grad} and \ref{p:grad1}.}
In this subsection, we will combine some results in \cite{KR16, KS15, KS17} to prove Proposition \ref{p:grad}.
Recall that we have assumed that $\nu:\R_+ \rightarrow \R_+$ is non-increasing differentiable  function satisfying $\int_{\R^d} \big(1 \land |y|^2 \big)\nu(|y|)dy<\infty$. 
In this subsection,
instead of the function $\Phi$, we mainly use 
\begin{equation}\label{d:varphi}
	 \varphi(r):= \begin{cases}\displaystyle \frac{r^2}{ \int_0^r s^{d+1}\nu(s)ds}, &\quad 0<r \le 1, \\
	\varphi(1)r^2, &\quad r > 1, \end{cases}
	\end{equation}Note that the integral $\int_0^r s^{d+1}\nu(s)ds$ above 
is finite because of our assumption  $\int_{\R^d} \big(1 \land |y|^2 \big)\nu(|y|)dy<\infty$.

To prove Propositions \ref{p:grad} and \ref{p:grad1} at once, we need to consider the following conditions on L\'evy measure $\nu(|y|)dy$ which is slightly more general than \eqref{c:nu}. 
We assume that there exist constants $a>0,0< \beta \le 1$ and $\ell \ge 0$ such that
\begin{equation}\label{c:nu1}
 \nu(r) \le ar^{-\ell} \exp(-br^\beta), \quad r>1.
\end{equation}
Also, we assume that there exist $a_3>0$ and $\alpha_3 \in (0,2]$ such that
\begin{equation}\label{e:varphi1}
a_3 \big(\frac Rr \big)^{\alpha_3} \le \frac{\varphi(R)}{\varphi(r)}, \quad 0<r \le R<\infty.
\end{equation}
For instance, when $X$ is an isotropic L\'evy process in Proposition \ref{p:grad} we have $\frac{s}{a\phi(s)} \le \nu(s)s^{d+1} \le \frac{as}{\phi(s)}$, which implies $\varphi(r) \asymp \Phi(r)$. Using this and Lemma \ref{e:Phi} we obtain \eqref{e:varphi1} with $\alpha_3=\alpha_1$. Thus, the conditions in Proposition \ref{p:grad} imply \eqref{c:nu1} and \eqref{e:varphi1}. 

Under \eqref{e:varphi1}, we have  $\varphi(r) \le cr^{\alpha_3}$ for $r \le 1$ so that 
$$c^{-1} r^{-\alpha_3} \le \int_0^r \frac{s^{d+1}}{r^2}\nu(s)ds \le \int_0^r s^{d-1}\nu(s)ds \le \int_0^1 s^{d-1}\nu(s)ds, \quad r \le 1. $$
Thus, letting $r \downarrow 0$ we obtain $\int_0^1 s^{d-1}\nu(s)ds=\infty.$
\noindent Now we record the counterpart of Lemma \ref{e:Phi}. Following the proof of Lemma \ref{e:Phi}, we obtain   
\begin{equation}\label{e:varphi2}
 \frac{\varphi(R)}{\varphi(r)} \le \big(\frac Rr \big)^{2}, \quad 0<r \le R.
\end{equation}
In addition, since $\nu$ is non-increasing, we have 
\begin{equation}\label{e:varphi}
\varphi(r)^{-1} ={r^{-2}\int_0^r s^{d+1}\nu(s)ds} \ge {r^{-2}\int_0^r s^{d+1} \nu(r)dr} = \frac{r^d\nu(r)}{d+2}, \quad r<1.
\end{equation}
In this subsection except the 
proofs of Propositions \ref{p:grad} and \ref{p:grad1}
we will always assume that $\nu$ satisfies \eqref{c:nu1} and \eqref{e:varphi1}. Let $X$ be the L\'evy process with L\'evy measure $\nu(|y|)dy$, and $\xi \mapsto \psi(|\xi|)$ be the characteristic exponent of $X$. First note that $\nu(\R^d)= \int_{\R^d} \nu(|y|)dy=\infty$ because $\int_0^1 s^{d-1}\nu(s)ds=\infty$. Also, since $X$ is isotropic, characteristic exponent of $X$ is also isotropic function. Define $\Psi(r):= \sup_{|y| \le r} \psi(|y|)$ and let $\mathcal{P}(r):= \int_{\R^d} \big( 1 \land \frac{|y|^2}{r^2} \big) \nu(|y|)dy$ be the Pruitt function for $X$. By \cite[Lemma 1 and Proposition 2]{BGR14}, we have that for $r>0$,
\begin{equation}\label{e:pru}
\frac{2}{\pi^2 d} \PP(r^{-1}) \le \psi(r) \le \Psi(r) \le \pi^2 \psi(r) \le 2\pi^2 \PP(r^{-1}), \quad r>0.
\end{equation}
Using \eqref{e:pru}, we can prove the following lemma. 
\begin{lemma}\label{l:usc-psi}
{Assume that $\nu(|y|)dy$ satisfies \eqref{c:nu1} and \eqref{e:varphi1}. Then, $\Psi(r)$ is comparable to $\varphi(r^{-1})^{-1}$, i.e., there exists a constant $c>0$ such that
\begin{equation}\label{e:Phi-psi}c^{-1}\varphi(r^{-1})^{-1} \le \Psi(r) \le c\varphi(r^{-1})^{-1}, \quad r>0.
\end{equation}}
\end{lemma}
\pf We claim that
\begin{equation}\label{e:pp}
\PP(r) \asymp \varphi(r)^{-1} \quad \mbox{for}  \quad r>0.
\end{equation}
 First assume $r \le 1$ and observe that
\begin{align*}
\PP(r) &= \int_{\R^d} \big(1 \land \frac{|z|^2}{r^2} \big) \nu(z) dz \\
&= c(d) \left( r^{-2} \int_0^r  s^{d+1} \nu(s)ds + \int_r^1 s^{d-1} \nu(s) ds + \int_1^\infty s^{d-1} \nu(s) ds \right) \\
&=: c(d) \big(I_{1} + I_{2} + I_{3} \big).
\end{align*}
By the definition of $\varphi$ we have
$ I_{1} ={\varphi(r)}^{-1}. $
{To estimate $I_2$, let us define $k:= \lfloor \frac{\log r}{\log 2} \rfloor$,
the largest integer smaller than or equal to $\frac{\log r}{\log 2}$. Then we have
\begin{align*}
 0 &\le I_{2} \le \sum_{i=0}^k \int_{2^i r}^{2^{i+1}r} s^{d-1}\nu(s)ds=: \sum_{i=0}^k I_{2i}.
 \end{align*}
Using \eqref{e:varphi1}, we have
\begin{align*}
 I_{2i} &\le (2^i r)^{-2} \int_{2^i r}^{2^{i+1}r} s^{d+1} \nu(s)ds \le  (2^ir)^{-2} \int_0^{2^{i+1}r} s^{d+1} \nu(s)ds = 4\varphi(2^{i+1} r)^{-1}  \le a_3 2^{2-\alpha_3(i+1)} \varphi(r)^{-1}.
 \end{align*}
 Thus,
 \begin{equation}\label{e:I2}
 0 \le I_2 \le \sum_{i=0}^k I_{2i} \le \frac{2^{2-\alpha_3}}{\varphi(r)} \sum_{i=0}^k 2^{-\alpha_3 i} \le \frac{c_1}{\varphi(r)}.
 \end{equation}}
Also, using \eqref{c:nu1} and \eqref{e:varphi1} we obtain 
$$ 0 \le I_{3} \le a \int_1^\infty s^{d-\ell-1} \exp(-bs^\beta) ds = c_2 \le \frac{ c_2 \varphi(1)}{a_3\varphi(r)}, $$
where we used $a_3 \le a_3 \big( \frac 1r \big)^{\alpha_3} \le \frac{\varphi(1)}{\varphi(r)}$ for the last inequality. Combining estimates of $I_1,I_2$ and $I_3$ we have proved the claim  \eqref{e:pp} for $r \le 1$. \\

\noindent Now assume $r > 1$. Then we have
\begin{align*}
\PP(r) &= \int_{\R^d} \big( 1 \land \frac{|z|^2}{r^2} \big) \nu(z)dz \\
&= c(d) \left( r^{-2} \int_0^1 s^{d+1} \nu(s)ds +  \int_1^\infty \big(1 \land \frac{s^2}{r^2} \big) s^{d-1}\nu(s)ds \right) \\
&:= c(d)(\varphi(r)^{-1} + I_4).
\end{align*}
Also, using \eqref{c:nu1} we have
$$ 0 \le I_4 \le \int_1^\infty \frac{s^2}{r^2} s^{d-1} \nu(s)ds \le ar^{-2} \int_1^\infty s^{d-\ell+1}\exp(-bs^\beta)ds \le c_3r^{-2}. $$
Using $\varphi(r)= \varphi(1)r^2$ for $r \ge 1$ we obtain that
$ \PP(r) \asymp r^{-2} \asymp \varphi(r)^{-1}$ for  $r > 1,$
which implies \eqref{e:pp} for $r >1 $. Therefore, \eqref{e:pp} holds for any $r>0$. Combining \eqref{e:pp} and \eqref{e:pru} we conclude the lemma. \qed

Using \eqref{e:Phi-psi}, \eqref{e:varphi1} and \eqref{e:varphi2} we obtain the following weak scaling condition for $\Psi$:  there exists a constant $c > 0$ such that
\begin{equation}\label{e:wsc-Psi}
c^{-1} \big( \frac{R}{r}\big)^{\alpha_3} \le \frac{\Psi(R)}{\Psi(r)} \le c \big( \frac{R}{r}\big)^2, \quad 0 < r \le R < \infty.
\end{equation}

Let $p_t(x)$ be a transition density of $X$. Since $X$ is isotropic, $x \mapsto p_t(x)$ is also isotropic function for any $t>0$. By an abuse of notation we also denote the radial part of the heat kernel $p_t(x)$  of $X$ as $p_t(r)$, $r>0$. 

To obtain gradient estimate for $p_t(x)$, we first follow the proof of \cite[Proposition 3.1]{KR16} to construct a $(d+2)$-dimensional L\'evy process $Y$ whose heat kernel estimate implies gradient estimate of $X$. \begin{lemma}\label{l:d+2}
	Assume that isotropic unimodal L\'evy measure $\nu$ satisfies \eqref{c:nu1} and \eqref{e:varphi1}. Then there exists an isotropic L\'evy process $Y$ in $\R^{d+2}$ such that its characteristic exponent is $\xi \mapsto \psi(|\xi|)$, $\xi \in \R^{d+2}$. Let $\nu_1(|x|)$ and $q_t(|x|)$ be the jumping kernel and heat kernel of $Y$, respectively. Then  for any $r>0$,
	\begin{equation}\label{e:p-d+2}
	q_t(r)= -\frac{1}{2 \pi r} \frac{d}{dr}p_t(r)
	\end{equation}
	and
	\begin{equation}\label{e:nu-d+2}
	\nu_1(r)=-\frac{1}{2\pi r} \nu'(r).
	\end{equation}
\end{lemma}
\pf The existence of $Y$ and \eqref{e:p-d+2} are immediately followed by \cite[Proposition 3.1]{KR16}. 
Note that using \eqref{e:pru} and \eqref{e:wsc-Psi} we have
$$ \lim_{\rho \rightarrow \infty} \frac{\psi(\rho)}{\log \rho}  \ge
\lim_{\rho \rightarrow \infty} \frac{\Psi(\rho) }{  \pi^2 \log \rho}  \ge \lim_{\rho \rightarrow \infty}\frac{c_1 \rho^{\alpha_3}}{\log \rho} = \infty,$$
which is one of the conditions in \cite[Proposition 3.1]{KR16}. For \eqref{e:nu-d+2}, we just need to follow the corresponding part in  the proof of \cite[Theorem 1.5]{KR16}. Here we provide a brief sketch for the proof for reader's convenience;
As in the proof of \cite[Theorem 1.5]{KR16}, without using
the assumption that $-\nu'(r)/r$ is non-increasing,  one can show that  
 there exists an  isotropic L\'evy process $X^{(d+2)}$ in $\R^{d+2}$ with jumping kernel $\nu_1(dy)$ and that the characteristic exponent of $X^{(d+2)}$ is $\psi(r)$. Thus, $X^{(d+2)}$ and $Y$ are identical in law, which concludes the proof. 
 To show this, only \cite[(8) and (9)]{KR16} are used, which follow directly from the fact that $\nu$ is isotropic, unimodal measure satisfying
$ \int_{\R^d} \big( |y|^2 \land 1) \nu(dy)< \infty. $
\qed

We emphasize here that we don't impose the condition \eqref{e:J2} on $\nu$. Thus  the function $r \to \nu_1(r)$ in the above lemma may not be non-increasing. 

{Now we are going to establish heat kernel estimates for the process  $Y$ obtained in Lemma \ref{l:d+2}, which will imply heat kernel estimate and gradient estimate of $X$ as a consequence of \eqref{e:p-d+2}. To do this,
we will  check conditions \textbf{(E), (D), (P)} and \textbf{(C)} (when $\beta < 1$) in \cite{KS17} for the process $X$ and $Y$, and apply \cite[Theorem 4]{KS17} and \cite[Theorem 1]{KS15}. }

First, we verify the condition \textbf{(E)} in \cite{KS17}. Recall $\Psi(r)= \sup_{|y| \le r} \psi(|y|)$.

\begin{lemma}\label{l:E}
Assume that isotropic unimodal L\'evy measure $\nu$ satisfies \eqref{c:nu1} and \eqref{e:varphi1}. Then for any $n,m \in \N$, 
there exists a constant $c=c(n,m)>0$ such that
$$ {\int_{\R^{n}} e^{-t \psi(|z|)} |z|^m dz \le c \Psi^{-1}(t^{-1})^{n+m}, \quad t>0.}  $$
\end{lemma}
\pf By \eqref{e:pru} and \eqref{e:wsc-Psi} we have that for $0<t $,
\begin{align*}
\int_{\R^n} e^{-t \psi(|z|)} |z|^m dz 
&\le c_1 \int_0^{\Psi^{-1}(t^{-1})} r^{n+m-1} dr +c_1\int_{\Psi^{-1}(t^{-1})}^\infty e^{-\pi^{-2} t \Psi(r)} r^{n+m-1}dr \\
&\le c_2  \Psi^{-1}(t^{-1})^{n+m} + c_1 \int_{\Psi^{-1}(t^{-1})}^\infty e^{-c_3 t \Psi(\Psi^{-1}(t^{-1})) (r/\Psi^{-1}(t^{-1}))^{\alpha_3}} r^{n+m-1} dr \\
&= \left(c_2 + c_1 \int_1^\infty e^{-c_3 s^{\alpha_1}}s^{n+m-1} dr \right)  \Psi^{-1}(t^{-1})^{n+m} = c_4  \Psi^{-1}(t^{-1})^{n+m},
\end{align*}
where we have used the change of variables with $s = \frac{r}{\Psi^{-1}(t^{-1})}$ in  the last line. \qed

{Note that Lemma \ref{l:E} for $(n,m)=(d,1)$ and $(n,m)=(d+2,1)$ implies the condition \textbf{(E)}  in \cite{KS17} for the process $X$ and $Y$, respectively.}

 For $0<\beta \le 1$ and $\ell \ge 0$,  we define non-increasing functions $f$ and $\tf$ by
 \begin{align}\label{d:f}
  f(r):= \begin{cases} \displaystyle \frac{\varphi(1)}{r^{d+1} \varphi(r)}, &\quad r \le 1, \\
 r^{-\ell-1}\exp(-br^\beta), &\quad r > 1 \end{cases} \quad \mbox{and} \quad \tf(r) := \begin{cases} \displaystyle \frac{\varphi(1)}{r^d \varphi(r)}, &\quad r \le 1, \\
 r^{-\ell}\exp(-br^\beta), &\quad r > 1 \end{cases}
 \end{align}
The functions $f$ and $\tf$ above are non-increasing since for any $0<r \le R \le 1$,
$$ \frac{1}{r^{d}\varphi(r)}= r^{-1}\int_0^r (\frac{s}{r})^{d+1}\nu(s)ds = \int_0^1 t^{d+1} \nu(rt)dt \ge  \int_0^1 t^{d+1} \nu(Rt)dt = \frac{1}{R^{d} \varphi(R)}. $$
Here we used that $\nu$ is nonincreasing. Note that by \eqref{c:nu1} and \eqref{e:varphi},
\begin{equation}\label{e:nu}
\frac{\nu(r)}{r} \le  cf(r) \quad \mbox{and} \quad \nu(r) \le c\tf(r) \qquad \mbox{for} \quad r>0
\end{equation}

{In the next lemma we verify the condition \textbf{(D)} in \cite{KS17} for both $X$ and $Y$. 
In fact, 
we are going to verify \textbf{(D)} for $X$ with the above $\tf$ and $\gamma=d$, while we use $f$ and $\gamma=d+1$ to verify \textbf{(D)} for $Y$. Let $B_d(x,r):=\{y \in \R^{d}: |x-y|<r\}$ and recall that  {\rm diam}$(A) = \sup \{ |x-y| : x,y \in A \}$ and $\nu_1(r)=-\frac{1}{2\pi r} \nu'(r)$.}

\begin{lemma}\label{l:D} 
Assume that $\nu$ satisfies \eqref{c:nu1} and \eqref{e:varphi1}. Then both $\nu(\R^{d})$ and
$\nu_1(\R^{d+2})=\int_{\R^{d+2}}\nu_1(|x|)dx$ are infinite, and there exists $c>0$ such that 
	\begin{equation}\label{e:(D)1}
\nu(A) \le c \tf(\delta(A)) [{\rm diam}(A)]^{d} , \quad A \in \BB(\R^{d}).
\end{equation}
and
\begin{equation}\label{e:(D)}
\nu_1(A) =\int_{A}\nu_1(|x|)dx  \le c f(\delta(A)) [{\rm diam}(A)]^{d+1} , \quad A \in \BB(\R^{d+2}). 
\end{equation}
for some $c>0$, where $\delta(A) := \inf\{ |y| : y \in A  \}$.  
\end{lemma}
\pf We have already showed that $\nu(B_d(0,1))=\nu(\R^d)=\infty$. For any $A \in \BB(\R^d)$, using \eqref{e:nu} we have
\begin{align*}
\nu(A) = \int_{A} \nu(|y|)dy \le  \nu(\delta(A))[{\rm diam}(A)]^{d} \le c\tf(\delta(A))[{\rm diam}(A)]^{d}.
\end{align*}
This concludes \eqref{e:(D)1}. \\

Using $\nu'(r) \le 0$, \eqref{c:nu}, the integration by parts and the fact $\nu(B_d(0,1))=\infty$ we have
\begin{align*}
\nu_1(\R^{d+2}) &\ge \int_{B_{d+2}(0,1)} \nu_1(|y|)dy =c(d) \int_0^1 r^{d+1} \nu_1(r)dr= c_1 \liminf_{\eps \downarrow 0} \int_\eps^1 -r^d \nu'(r)dr \\
&= c_1 \liminf_{\eps \downarrow 0} \left( -[r^d \nu(r)]_\eps^1 + d \int_\eps^1 r^{d-1} \nu(r)dr \right) \\ &= c_1 \liminf_{\eps \downarrow 0} \big( \eps^d \nu(\eps) +d \int_\eps^1 r^{d-1} \nu(r)dr - \nu(1) \big) \ge c_2 \nu(B_d(0,1))- c_1\nu(1) = \infty.
\end{align*}
Now it remains to prove \eqref{e:(D)}.  First observe that using the integration by parts, we have that for any $0<r<R$,
\begin{align}
\begin{split}\label{e:(1).1}
\int_r^R s^{d+1}\nu_1(s)ds &= -\frac{1}{2 \pi} \int_r^R s^{d} \nu'(s)ds = \frac{1}{2 \pi} \left(-[s^d \nu(s)]_r^R + d\int_r^R s^{d-1}\nu(s)ds \right)\\ &\le \frac{1}{2 \pi} \left( r^d \nu(r) + \nu(r) \, d \int_r^R s^{d-1}ds \right) = \frac{1}{2 \pi} \nu(r) R^d
\end{split}
\end{align}
where we used that $\nu$ is non-increasing. Now denote $r:= \delta(A)$ and $ l:=$diam$(A)$.

 When $l \ge r/2$, using $A \subset \{ y \in \R^{d+2}: r \le |y| \le r+l  \}$ we obtain
\begin{align*}
\nu_1(A) &\le \nu_1 ( \{ y: r \le |y| \le r+l  \} ) = c(d) \int_r^{r+l} s^{d+1} \nu_1(s)ds \\
&\le  \frac{c(d) }{2\pi} \nu(r)(r+l)^{d} \le c_1  \frac{\nu(r)}{r} l^{d+l} \le c_3 f(r) l^{d+1},
\end{align*}
where we used \eqref{e:(1).1} and \eqref{e:nu} for the last line. 

 When $l < r/2$, choose a point $y_0 \in \bar{A}$ with $|y_0|=r$. Since $A \subset B_{d+2}(y_0,l) \backslash B_{d+2}(0,r)$, there exists $c_4=c_4(d)>0$ such that 
$$ \int_{|y|=s} \1_A(y)  \sigma(dy) \le c_4 l^{d+1}$$
for any $s \in [r,r+l]$. Thus, by \eqref{e:(1).1} and \eqref{e:nu} we have
\begin{align*}
\nu_1(A) &\le \nu_1 ( B(y_0,l) \backslash B(0,r))\le  c_5\int_r^{r+l} l^{d+1} \nu_1(s)ds \le c_5 \frac{l^{d+1}}{r^{d+1}} \int_r^{r+l} s^{d+1}\nu_1(s)ds \\
&\le \frac{c_5}{2 \pi} \frac{l^{d+1}}{r^{d+1}} (r+l)^d \nu(r) \le c_6 l^{d+1} \frac{\nu(r)}{r} \le c_7 f(r) l^{d+1},
\end{align*}
which proves \eqref{e:(D)}. 
\qed

Recall $\Psi(r)= \sup_{|y| \le r} \psi(|y|)$. 
\begin{lemma}\label{l:P0}
Assume that $\nu$ satisfies \eqref{c:nu1} and \eqref{e:varphi1}.
	For every $\kappa<1$, there exists $c=c(\kappa)>0$ such that 
	\begin{equation}\label{e:tnu1}
	\int_{\{y\in\R^{d}:|y|>r\}} \exp \big(b\kappa |y|^\beta\big) \nu(dy) \le c\Psi( \frac{1}{r}), \quad r>0
	\end{equation}
		and
		\begin{equation}\label{e:nu1}
	\int_{\{y\in\R^{d+2}:|y|>r\}} \exp \big(b\kappa |y|^\beta\big) \nu_1(dy) \le c\Psi( \frac{1}{r}), \quad r>0
	\end{equation}

\end{lemma}
\pf 
Since \eqref{e:tnu1} can be derived directly from the estimate of $I_2$ below, 
we only prove \eqref{e:nu1} here. Using the  integration by parts, we have
\begin{align*}
\int_{|y|>r}  \exp(b\kappa |y|^\beta) \nu_1(dy) &= c(d) \int_r^\infty \exp(b\kappa t^\beta) t^d (-\nu'(t))dt \\
&= c(d)\left( \left[\exp(b\kappa t^\beta) t^d (-\nu(t))\right]_r^\infty + \int_r^\infty (\exp(b\kappa t^\beta)t^d)' \nu(t)dt \right). \\
&:= c(d) \big( I_1 + I_2 \big).
\end{align*}
For $I_1$, by \eqref{e:nu} $\displaystyle \lim_{t \rightarrow \infty} \exp(b\kappa t^\beta)t^d \nu(t) \le \lim_{t \rightarrow \infty} a\exp(-b(1-\kappa) t^\beta)) t^{d-\ell} =0 $
so
$ I_1 = \exp(b\kappa r^\beta)r^d \nu(r)$ $\le c_1  \varphi(r)^{-1}.$
Now let us estimate $I_2$. First we observe that
$$\frac{d}{dt} \big( \exp(b\kappa t^\beta)t^d \big) \le c_2 \begin{cases} t^{d-1}, &\quad t \le 1 \\
\exp(b\kappa t^\beta)t^{d+\beta-1}, &\quad t > 1.\end{cases} $$
Thus, for $r \ge 1$ we have
\begin{align*}
\int_r^\infty (\exp(b\kappa t^\beta)t^d)' \nu(t) dt \le c_2 \int_r^\infty \exp(-b(1-\kappa)t^\beta ) t^{d-\ell+\beta-1} dt \le c_3 r^{-2}= \frac{c_3 \varphi(1)}{\varphi(r)} .
\end{align*}
For $r < 1$, using above estimate, \eqref{e:I2} and \eqref{e:varphi1} we get
\begin{align*}
\int_r^\infty (\exp(b\kappa t^\beta)t^d)' \nu(t) dt &= \left(\int_r^{1} + \int_{1}^\infty  \right) (\exp(b\kappa t^\beta)t^d)' \nu(t) dt \\
&\le c_2 \left(\int_r^1 t^{d-1}\nu(t)dt + \int_1^\infty \exp(-b(1-\kappa)t^\beta)t^{d-\ell+\beta-1} dt \right) \\
&\le \frac{c_4}{\varphi(r)} + c_3 \le \frac{c_5}{\varphi(r)}
\end{align*}
Combining above two inequalities and \eqref{e:Phi-psi}, we obtain
$ I_1 + I_2 \le c_6 \Psi(\frac{1}{r})$. Therefore, we have proved the lemma.
\qed

Using Lemma \ref{l:P0}, we verify the condition \textbf{(P)} in \cite{KS17} for both $X$ and $Y$. We continue to use the non-increasing functions $f$ and $\tf$ defined in \eqref{d:f}.
\begin{lemma}\label{l:P}
Assume that isotropic unimodal L\'evy measure $\nu$ satisfies \eqref{c:nu1} and \eqref{e:varphi1}. 
	Then, there exists $c>0$ such that
	\begin{equation}\label{e:P1}
	\int_{\{y \in \R^{d}:|y|>r\}} \tf \left( s \lor |y| - \frac{|y|}{2} \right) \nu(dy) \le c\tf(s) \Psi(\frac{1}{r}), \quad r,s>0
	\end{equation}
	and
		\begin{equation}\label{e:P}
	\int_{\{y \in \R^{d+2}:|y|>r\}} f \left( s \lor |y| - \frac{|y|}{2} \right) \nu_1(dy) \le cf(s) \Psi(\frac{1}{r}), \quad r,s>0
	\end{equation}
\end{lemma}
\pf We only prove \eqref{e:P} here, since \eqref{e:P1} can be verified  similarly. We claim that for any $0<\beta \le 1$, there exists $c_1>0$ such that for any $s,t>0$, 
\begin{equation}\label{e:f}
f(s \lor t - \frac{t}{2}) \le c_1 f(s) \exp(b\kappa t^\beta)
\end{equation}
where $\kappa = \frac{1}{2}(2^{-\beta}+1)$. First we define
$$  f_1(r):= \begin{cases}\frac{\varphi(1)}{r^{d+1}\varphi(r)}, &\quad r \le 2 \\ r^{-\ell-1}\exp(-br^\beta), &\quad r>2. \end{cases}$$
Then, since $f(r)=f_1(r)$ for $r \in (0,1] \cup (2,\infty)$ we have
\begin{equation}\label{e:f_1}   c_2^{-1} f(r) \le f_1(r) \le c_2f(r), \quad r>0.
\end{equation}
Now assume $s \lor t > 2$. Then, using $1 \lor \frac{s}{2} \le s \lor t- \frac{t}{2}$ and triangular inequality,
\begin{align*}
f(s \lor t - \frac{t}{2})  &= (s \lor t - \frac{t}{2})^{-\ell-1} \exp\big(-b (s \lor t - \frac{t}{2})^\beta   \big) \le (1 \lor \frac{s}{2})^{-\ell-1} \exp(-b(s-\frac{t}{2})^\beta) \\ &\le   (1 \lor \frac{s}{2})^{-\ell-1} \exp(-bs^\beta) \exp(b(\frac{t}{2})^\beta) \le c_3f(s)\exp(b(\frac{t}{2})^\beta).
\end{align*}
Here in the last inequality we used $\ell \ge 0$ and $\exp(-bs^\beta) \le c_3f(s)$ for $0<s \le 2$. 
When $s \le 2$ and $t \le 2$, using \eqref{e:f_1}, \eqref{e:varphi1} and \eqref{e:varphi2} with $s \lor t - \frac{t}{2} \ge \frac{s}{2}$ we obtain
\begin{align*}
f(s \lor t - \frac{t}{2}) \le c_2 f_1(s \lor t - \frac{t}{2}) = \frac{c_2\varphi(1)}{(s \lor t - \frac{t}{2})^{d+1}\varphi(s \lor t - \frac{t}{2})} \le \frac{c_4 \varphi(1)}{s^{d+1}\varphi(s)} \le c_4 f_1(s) \le c_5f(s). 
\end{align*}
Here we used $\frac{\varphi(s)}{\varphi(s \lor t - \frac{t}{2})} = \frac{\varphi(s)}{\varphi(s/2)} \frac{\varphi(s/2)}{\varphi(s \lor t - \frac{t}{2})} \le 4a_3^{-1}$ which follows from \eqref{e:varphi1} and \eqref{e:varphi2}. Thus, we conclude \eqref{e:f}. Combining  \eqref{e:f} and Lemma \ref{l:P0}, we have proved the lemma. \qed

Now we obtain a priori heat kernel estimates for the process $X$ and $Y$. To state the results, we need to define \textit{generalized inverse of $\varphi$} by
$ \vinv(t) := \inf\{ s >0 : \varphi(s) \ge t\}.$
Using \eqref{e:varphi1} and \cite[Remark 4]{BGR14}, we  obtain
\begin{equation}\label{e:wsc-vinv}
c^{-1}\big(\frac{R}{r}\big)^{1/2}  \le \frac{\vinv(R)}{\vinv(r)} \le c \big(\frac{R}{r}\big)^{1/\alpha_3}
\end{equation}
and
\begin{equation}\label{e:vinv}
a_3^{-1} \varphi(\vinv(r)) \le r \le a_3 \varphi(\vinv(r)),
\end{equation}
which are counterparts of \eqref{e:wsc-inv}. First we apply \cite[Theorem 3]{KS15} to obtain the regularity of the transition density $p_t(x)$ of $X$.

\begin{prop}\label{p:reg}
	Let $X$ be an isotropic unimodal L\'evy process in $\R^d$ with jumping kernel $\nu(|y|)dy$ satisfying \eqref{c:nu1} and \eqref{e:varphi1} with $0<\beta \le 1$. Then $x \to p_t(x) \in C^\infty_b(\R^d)$ and for any $k \in \N_0$ there exists $c_k>0$ such that
	\begin{equation}\label{e:reg}
	|\nabla^k_x p_t(x)| \le c_k  \vinv(t)^{-k} \left( \vinv(t)^{-d} \land \frac{t}{|x|^d\varphi(|x|)}\right)
	\end{equation}
	for any $t>0$ and $x \in \R^d$.
\end{prop}
\pf 
Define $h(t):= \frac{1}{\Psi^{-1}(t^{-1})}$ as in \cite{KS15}. Note that by \eqref{e:Phi-psi} and \eqref{e:vinv} we have
\begin{equation}\label{e:h}
h(t) \asymp \vinv(t), \quad t>0.
 \end{equation}
{Applying \cite[Theorem 3]{KS15} for the process $X$}, $p_t(x) \in C^\infty_b(\R^d)$ and for any $k \in \N$, $\gamma \in [1,d]$ and $n>\gamma$ we have constants $c_{k,n}$ satisfying
$$ |\nabla^k_x p_t(x)| \le c_{k,n} (h(t))^{-d-k} \min \left\{ 1, \frac{t [h(t)]^\gamma}{|x|^{\gamma}  \varphi(|x|)}e^{-b(|x|/4)^\beta} + \left( 1+ \frac{|x|}{h(t)} \right)^{-n} \right\}. $$
Note that we already verified \cite[(8)]{KS15} at Lemma \ref{l:E}.  Thus, using $h(t) \asymp \vinv(t)$ we obtain
$$ |\nabla^k_x p_t(x)| \le \tilde{c}_{k,n} \vinv(t)^{-d-k} $$
Also, taking $\gamma=d$, $n=d+2$ and using $h(t) \asymp \vinv(t)$ we get
\begin{align*}
|\nabla^k_x p_t(x)| &\le c_{k,n} \left( h(t)^{-k} \, \frac{t}{|x|^d \varphi(|x|)}  e^{-b(|x|/4)^\beta} + h(t)^{-k} |x|^{d} \left( 1+ \frac{|x|}{h(t)} \right)^{-2} \right)  \\
&\le c \vinv(t)^{-k} \left( \frac{t}{|x|^d \varphi(|x|)} +   \big(\frac{\vinv(t)}{|x|} \land 1\big)^2 \, \frac{1}{|x|^d} \right) \\
&\le   \tilde{c}_{k,n} \vinv(t)^{-k} \, \frac{t}{|x|^d \varphi(|x|)}.
\end{align*}
The last inequality is straightforward when $|x| < \vinv(t)$ and it  follows from \eqref{e:wsc-vinv} and \eqref{e:vinv} when $|x| \ge \vinv(t)$.  Therefore, we conclude that
\begin{equation*}
|\nabla^k_x p_t(x)| \le c_k  \vinv(t)^{-k} \left( \vinv(t)^{-d} \land \frac{t}{|x|^d\varphi(|x|)}\right).
\end{equation*}
\qed

{Note that the gradient estimates in Proposition \ref{p:reg} is same as the ones in \cite[Proposition 3.2]{KSV} except that  the gradient estimates in \cite[Proposition 3.2]{KSV} is for $t \le T$ (see Remark \ref{r:KSV}). }

 Combining above estimates with Lemmas \ref{l:E}, \ref{l:D} and \ref{l:P}, we can apply \cite[Theorem 1]{KS15} for the process $X$ and $Y$. Here is the result.
\begin{lemma}\label{l:KS15}
Assume that $\nu$ satisfies \eqref{c:nu1} and \eqref{e:varphi1} and $\beta=1$. Then for any $T \ge 1$, there exists a constant $c>0$ such that
\begin{equation}\label{e:beta=1}
p_t(x) \le ct \exp(-\frac{b}{4}|x|) \quad \mbox{and} \quad q_t(x) \le ct\vinv(t)^{-1} \exp(-\frac{b}{4}|x|)
\end{equation}
for any $0<t \le T$ and $|x| > \vinv(T)$.
\end{lemma}
\pf Define $h(t):= \frac{1}{\Psi^{-1}(t^{-1})}$ as in \cite{KS15} and denote $q_t(|x|)=q_t(x)$. 
Applying  Lemmas \ref{l:E}, \ref{l:D} and \ref{l:P} to \cite[Thoerem 1]{KS15} for the process $Y$
 in Lemma \ref{l:d+2}, we have
\begin{align*}
 &q_t(r) \le c_1 h(t)^{-1} \left( h(t)^{-d-1} \land \left[tf(r/4)   + h(t)^{-d-1} \exp\big(-c_2 \frac{r}{h(t)} \log (1+\frac{r}{h(t)} )\big)\right] \right) \\
 &\le c_3 \vinv(t)^{-1} \left(  \vinv(t)^{-d-1} \land  \left[ tf(r/4) + \vinv(t)^{-d-1} \exp\big(-c_4  \frac{r}{\vinv(t)} \log(1+c_5\frac{r}{\vinv(t)}) \big) \right]  \right)
\end{align*}
for any $t, r>0$. 

First observe that using $ f(\frac{r}{4}) = (\frac{r}{4})^{-\ell-1}\exp(-\frac{b}{4}r)$ for $r>4$, we obtain \begin{equation}\label{e:h2}t\vinv(t)^{-1}f(\frac{r}{4}) \le c_6tr^{-\ell-1}\exp(-\frac{b}{4}r) \le c_7 t \exp(-\frac{b}{4}r), \quad r>4.
\end{equation}
 Let $c(T)>4$ be a constant which is large enough to satisfy $$\frac{c_4}{2\vinv(T)} \log(1+ c_5\frac{c(T)}{\vinv(T)}) \ge \frac{b}{4}, \quad r>1.$$
Then using \eqref{e:wsc-vinv} in the second inequality, for any $0<t \le T$ and $r>c(T)$ we have
 \begin{align*}
 &\vinv(t)^{-d-1} \exp\big(-c_4 \frac{r}{\vinv(t)} \log (1+c_5\frac{r}{\vinv(t)} )\big) 
 \nn \\ &\le \vinv(t)^{-d-1} \exp\big(-\frac{c_4r}{2\vinv(t)}  \log(1+ c_5\frac{c(T)}{\vinv(T)})\big) \exp\big(-\frac{c_4r}{2\vinv(t)}  \log(1+ c_5\frac{c(T)}{\vinv(T)})\big)  
\nn  \\ &\le \vinv(t)^{-d-1} \exp\big(-c_8 \frac{r}{\vinv(t)} \big)\exp\big(-\frac{c_4r}{2\vinv(T)}  \log(1+ c_5\frac{c(T)}{\vinv(T)})\big)  \\
 &\le \vinv(t)^{-d-1} \exp(-c_8 \frac{r}{\vinv(t)}-\frac{b}{4}r) \le \frac{c_9t}{r^{d+1}\varphi(r)} \exp(-\frac{b}{4}r) \le c_9t\exp(-\frac{b}{4}r) \nn
 \end{align*}
 where $\displaystyle c_9 = \sup_{s \ge 1} s^{d+1} \exp(-c_8 s)< \infty$.
 Thus, 
 \begin{equation*}
 q_t(r) \le c_{10} t\vinv(t)^{-1} \exp(-\frac{b}{4}r), \quad r>c(T)= \varphi(\vinv(c(T))).
 \end{equation*}
 Meanwhile, by \eqref{e:p-d+2} and \eqref{e:reg} we have
 $$ q_t(r)= \frac{1}{2\pi r} |\frac{r}{dr} p_t(r)| \le \frac{c_{11} t\vinv(t)^{-1}}{r^{d+1} \varphi(r)} \le c_{12} t\vinv(t)^{-1} \exp(-\frac{b}{4}r+ \frac{bc(T)}{4})$$
 for $\vinv(T) < r \le c(T)$.
Therefore, combining above two estimates we conclude the estimate on $q$ in \eqref{e:beta=1}. 

Note that,  applying  Lemmas \ref{l:E}, \ref{l:D} and \ref{l:P} to \cite[Thoerem 1]{KS15} for the process $X$ and using $h(t) \asymp \vinv(t)$ we have
\begin{equation}\label{e:11}
p_t(r) \le c_{10} \left( \vinv(t)^{-d} \land \left[ t \tf(r/4) + \vinv(t)^{-d} \exp(-c_{11} \frac{r}{\vinv(t)} \log(1+ c_{12} \frac{r}{\vinv(t)}) \right] \right).
\end{equation}
for any $t,r>0$. Using \eqref{e:11}, the estimate on $p$ in \eqref{e:beta=1} can be verified similarly. \qed

Now we check condition \textbf{(C)} in \cite{KS17} with $r_0=1$, $t_p=\infty$ and $\gamma=d$ for $X$ ($\gamma=d+1$ for $Y$, respectively). We need additional condition $0<\beta<1$ to verify the condition \textbf{(C)}.

\begin{lemma}\label{l:C}
Assume $\nu$ satisfies \eqref{c:nu1} and \eqref{e:varphi1} with $0<\beta<1$. Then, there exists constant $c>0$ such that for every $|x| \ge 2$ and $r \in (0,1]$,
\begin{align}
\label{e:C1}
\tf(r) \le cr^{-d}\Psi( \frac{1}{r}), \qquad \int_{\{y \in \R^d:|x-y| \ge 1, |y| > r\}} \tf(|x-y|) \nu(dy) \le c \Psi(\frac{1}{r}) \tf(|x|) 
\end{align}
and
\begin{align}
\label{e:C}
f(r) \le cr^{-d-1}\Psi( \frac{1}{r}), \qquad \int_{\{y \in \R^{d+2} : |x-y| \ge 1, |y| > r\}} f(|x-y|) \nu_1(dy) \le c \Psi(\frac{1}{r}) f(|x|). 
\end{align}
\end{lemma}
\pf 
The first inequalities  in \eqref{e:C1} and \eqref{e:C} immediately follow from \eqref{e:Phi-psi} and \eqref{d:f}.

Let us show the second inequality  in \eqref{e:C}. When $|x-y| \ge \frac{|x|}{2}$, using \eqref{e:exp2} and triangular inequality,
we have $|x|^\beta \le |x-y|^\beta + (2^\beta-1)|y|^\beta. $
Thus, using this inequality and Lemma \ref{l:P0} we obtain
\begin{align*}
\int_{|x-y| \ge \frac{|x|}{2}, |y| > r} &f(|x-y|) \nu_1(dy) = \int_{|x-y| \ge \frac{|x|}{2}, |y| > r} |x-y|^{-\ell-1}\exp(-b|x-y|^\beta) \nu_1(dy) \\
&\le \big(\frac{|x|}{2}\big)^{-\ell-1} \int_{|y|>r} \exp(-b|x|^\beta) \exp(b(2^\beta-1)|y|^\beta)  \nu_1(dy) \\
&=f(|x|) \int_{|y|>r} \exp(b(2^\beta-1)|y|^\beta)  \nu_1(dy) \le c_1f(|x|) \Psi(\frac 1r).
\end{align*}  
 So, it suffices to show that there exists a constant $c_2>0$ such that for every $|x| \ge 2$,
\begin{align}
\label{e:claim1}
\int_{1 \le |x-y| \le \frac{|x|}{2}} f(|x-y|) \nu_1(dy) \le c_2  f(|x|).
\end{align}
To show this, we will divide the set $D:=\{ y: 1 \le |x-y| \le \frac{|x|}{2} \}$ into cubes with diameter 1. Let $x=(x_1,...,x_{d+2})$. For $(a_1,...,a_{d+2}) \in \Z^{d+2}$, we define $a:=(\sqrt{d+2})^{-1}(a_1,...,a_{d+2})$, and let
 $$  C_{a} := \prod_{i=1}^{d+2}[x_i+\frac{2a_i-1}{2\sqrt {d+2}},x_i+\frac{2a_i+1}{2\sqrt {d+2}})$$
be a cube with length $(\sqrt{d+2})^{-1}$. Since diam$(C_a) = 1$ and $x + a$ is the center of cube $C_a$, for any $|a| \le \frac{|x|+1}{2}$ we have $c_5>0$ independent of $a$ such that
\begin{align*}
\nu_1(C_a \cap D) &\le c_3 f(\delta(C_a \cap D)) \le c_3 f\left(\big(|x|-|a|- \frac{1}{2}\big) \lor \frac{|x|}{2}\right) \\ &\le c_4\big(|x|-|a|)^{-\ell-1}\exp\left(-b\left||x|-|a|\right|^\beta\right) \le c_5 |x|^{-\ell-1}\exp\left(-b \big(|x|-|a|\big)^\beta\right)
\end{align*}
where we used Lemma \ref{l:D} for the first inequality and triangular inequality for the second line. Thus, using $|a|-  \frac{1}{2} \le |x-y|$ on $C_a$ and $$ D \subset  \displaystyle \bigcup_{1 \le |a| \le \frac{|x|+1}{2}} C_a,  $$ we obtain
\begin{align*}
\int_{1 \le |x-y| \le \frac{|x|}{2}} f(|x-y|) \nu_1(dy) 
&\le \sum_{1\le |a| \le \frac{|x|+1}{2}}  \int_{C_a \cap D} |x-y|^{-\ell-1}\exp(-b|x-y|^\beta)\nu_1(dy) \\
&\le \sum_{1\le |a| \le \frac{|x|+1}{2}} (|a|-\frac{1}{2})^{-\ell-1}\exp(-b(|a|- \frac{1}{2})^\beta)\nu_1(C_a \cap D) \\
&\le c_6 |x|^{-\ell-1} \sum_{1\le |a| \le \frac{|x|}{2}+1} |a|^{-\ell-1}\exp(-b|a|^\beta) \exp(-b(|x|-|a|)^\beta).
\end{align*}
Since $|a| \le \frac{|x|+1}{2}$, by \eqref{e:exp2} we have
$|a|^\beta + (|x|-|a|)^\beta +1 \ge |a|^\beta + (|x|+1-|a|)^\beta \ge |x|^\beta + (2-2^\beta)|a|^\beta. $
Thus,
\begin{align*}
&|x|^{-\ell-1}\sum_{1\le |a| \le \frac{|x|+1}{2}} \exp(-b|a|^\beta) \exp(-b(|x|-|a|)^\beta) \\ &\le c_7 |x|^{-\ell-1}\exp(-b|x|^\beta) \sum_{a \in \Z^d \setminus \{0\}} |a|^{-\ell-1} \exp(-b(2-2^\beta)|a|^\beta) \le c_8 f(|x|). 
\end{align*}
Combining above inequalities and using \eqref{e:Phi-psi}, we arrive \eqref{e:claim1}.
Therefore, we conclude that the second inequality in \eqref{e:C} holds. 

The second inequality in \eqref{e:C1} can be verified similarly so skip the proof.
 \qed

\noindent Now we have that conditions \textbf{(E), (D)} and \textbf{(C)} in \cite{KS17} holds for  the process $Y$ when $\nu$ satisfies \eqref{c:nu1} and \eqref{e:varphi1} with $0<\beta<1$. Thus, we can apply \cite[Thereoem 4]{KS17} for both $X$ and $Y$. 
\begin{lemma}\label{l:KS17}
	Let $T \ge 1$ and assume that $\nu$ satisfies \eqref{c:nu1} and \eqref{e:varphi1} with $0<\beta<1$. Then, there exists a constant $c>0$ such that
	\begin{equation}\label{e:p-d}
	p_t(r) \le ct r^{-\ell} \exp(-br^\beta)
	\end{equation}
	and
	\begin{equation}\label{e:q-d+2}
\left|\frac{d}{dr}p_t(r)\right| \le ct \vinv(t)^{-1} r^{-\ell} \exp(-br^\beta)
	\end{equation}
	for any $0<t \le T$ and $r \ge 4$. 
\end{lemma}	 
\pf Applying \cite[Theorem 4]{KS17} for $Y$ and \eqref{e:Phi-psi} we have that for $0<t \le t_p= T$ and $r \ge 4r_0=4$,
$$ q_t(r)  \le c_1t\vinv(t)^{-1}f(r)= c_1t\vinv(t)^{-1}r^{-\ell-1}\exp(-br^\beta). $$
Combining this with \eqref{e:p-d+2}, 
$|\frac{d}{dr}p_t(r)| \le  2 \pi r q_t(r)   \le c_2 t\vinv(t)^{-1}r^{-\ell}\exp(-br^\beta). $ This concludes \eqref{e:q-d+2}. \eqref{e:p-d} immediately follows from applying \cite[Theorem 4]{KS17} for $X$.
 \qed

For reader's convenience, we put the heat kernel estimates and gradient estimates in Proposition \ref{p:reg}, and Lemmas \ref{l:KS15} and \ref{l:KS17} together into one proposition.  

\begin{prop}\label{p:hke}
Let $X$ be an isotropic unimodal L\'evy process in $\R^d$ with jumping kernel $\nu(|y|)dy$ satisfying \eqref{c:nu1} and \eqref{e:varphi1}. Then, $x \mapsto p_t(x) \in C_b^\infty (\R^d)$ and the following holds. \\
(a) There exists a constant $c_1>0$ such that
$$ |\nabla_x^k p_t(x)| \le c_1 \vinv(t)^{-k} \left( \vinv(t)^{-d} \land \frac{t}{r^d\varphi(r)}\right), \quad t>0, \, x\in \R^d \quad \mbox{and} \quad k \in \N_0. $$
(b) Assume $\beta=1$. Then for any $T \ge 1$, there exists a constant $c_2>0$ such that 
$$ |\nabla_x^k p_t(x)| \le c_2 t \vinv(t)^{-k} \exp(-\frac{b}{4}r), \quad t \in (0,T], \, |x| > \vinv(T) \quad \mbox{and} \quad k=0,1. $$
(c) Assume $0<\beta<1$. Then for any $T \ge 1$, there exists a constant $c_3>0$ such that 
$$ |\nabla_x^k p_t(x)| \le c_3 t \vinv(t)^{-k} r^{-\ell}\exp(-br^\beta), \quad t \in (0,T], \, |x| > \vinv(T) \quad \mbox{and} \quad k=0,1. $$
\end{prop}
\pf (a) and (b) immediately follow from Proposition \ref{p:reg} and Lemma \ref{l:KS15}, respectively. \\
(c) Observe that for any $t \in (0,T]$, $\vinv(T) <|x| \le 4$ and $k=0,1$,
\begin{align*}
|\nabla_x^k p_t(x)| \le c_1 \vinv(t)^{-k}  \frac{t}{r^d\varphi(r)} \le c_4 t \vinv(t)^{-k} r^{-\ell} \exp(-br^\beta).
\end{align*}
This and Lemma \ref{l:KS17} finish the proof. \qed

Now we are ready to prove Propositions \ref{p:grad} and  \ref{p:grad1}.

 \noindent \textbf{Proof of Proposition \ref{p:grad}.}  Now assume that $X$ is an isotropic L\'evy process in Proposition \ref{p:grad} with L\'evy measure $\nu(|y|)dy$. Recall that $\varphi(r) \asymp \Phi(r)$, and $\nu$ satisfies \eqref{c:nu1} with $\ell=0$ and \eqref{e:varphi1}. Therefore, we can apply results in Proposition \ref{p:hke} with ths function $\Phi$ instead of $\varphi$. Using Proposition \ref{p:hke} and \eqref{e:rho_T}, we conclude that for  any $t \in (0,T]$ and $x \in \R^d$ 
 $$ | \nabla_x^k p_t(x)| \le c_1 t \vinv(t)^{-k}\G_T(t,x) \le c_2 t \G_{-k}^0(t,x), \quad  k=0,1.$$ \qed

 \noindent \textbf{Proof of Proposition \ref{p:grad1}.} \eqref{e:grad2} for $k=0,1$ and that $t \mapsto p(t,x)$ is in $C_b^\infty(\R^d)$ immediately follow from Proposition \ref{p:grad}. 

Now it suffices to prove \eqref{e:grad2} when $k=2$. Let $X$ be an isotropic unimodal L\'evy process with jumping kernel $J(|x|)dx$ satisfying \eqref{e:J1} with $0<\beta \le 1$ and \eqref{e:J2}, and let $\psi(|x|)= \psi(x)$ be a characteristic exponent of $X$. By Lemma \ref{l:d+2}, there exists an isotropic L\'evy process $Y$ in $\R^{d+2}$ with characteristric exponent $\psi(r)$ satisfying \eqref{e:p-d+2} and \eqref{e:nu-d+2}. In particular, by \eqref{e:J2} and \eqref{e:nu-d+2}, $Y$ is unimodal. Denote $J_1(|x|)dx$ and $q_t(|x|)$ be the L\'evy density and transition density of $Y$ respectively. Using \eqref{e:nu-d+2} we have
$$
2\pi \int_s^r J_1(t)dt =- \int_s^r\frac{J'(t)}{t}dt = -[ \frac{J(t)}{t}]_s^r - \int_s^r  \frac{J(t)}{t^2}dt 
= \frac{J(s)}{s} - \frac{J(r)}{r} - \int_s^r \frac{J(t)}{t^2}dt.
$$ 
Since $J_1$ is non-increasing by \eqref{e:J2}, we obtain that for any $0<s \le r$, 
\begin{align}\label{e:3.2-1}
(r-s)J_1(r) &\le \int_s^r J_1(t)dt \le  \frac{J(s)}{2\pi s}
\end{align}
and
\begin{align}\label{e:3.2-2}
(r-s)J_1(s) 
\ge \int_s^r J_1(r)dr \ge \frac{1}{2\pi} \left( \frac{J(s)}{s} - \frac{J(r)}{r} - \int_s^r \frac{J(s)}{t^2}dt\right) = \frac{1}{2\pi r} (J(s)-J(r)) .
\end{align}

We claim that there exists a constant $c>0$ such that
\begin{align}\label{e:ned2}
{\frac{c^{-1}}{r^{d+2}\phi(r)} \le J_1(r) \le \frac{c}{r^{d+2}\phi(r/2)}, \quad r \le 1 \quad \mbox{and} \quad J_1(r) \le cr^{-1}\exp(-br^\beta), \quad r>1.}
\end{align}

For $r \le 1$, letting $s= \frac{r}{2}$ in \eqref{e:3.2-1} we have
$
J_1(r) \le \frac{J(r/2)}{\pi r^2} \le \frac{c_1}{r^{d+2}\phi(r/2)}
$ by \eqref{e:J1}.
Also, taking $(r,s)=(Cr,r)$ with constant $C = \big(\frac{2}{a^2}\big)^{2/d}>1$ in \eqref{e:3.2-2} we have
\begin{align*}
J_1(r) &\ge \frac{J(r)-J(Cr)}{2\pi C(C-1)r^2} \ge \frac{1}{2\pi C(C-1)r^2} \big( \frac{a}{r^d \phi(r)} - \frac{1}{a (Cr)^d \phi(Cr)} \big) \\ &\ge \frac{1}{C(C-1)} \left( a - \frac{1}{aC^d} \right) \frac{1}{r^{d+2} \phi(r)} = \frac{a}{2C(C-1)} \frac{1}{r^{d+2}\phi(r)},
\end{align*}
where we used $\phi(Cr) \ge \phi(r)$ and $a-\frac{1}{aC^d} = \frac{a}{2}$ in the second line. 

When $r>1$, letting $s=r-1$ in \eqref{e:3.2-1} we have
\begin{align*}
J_1(r) \le \frac{J(r-1)}{r} \le \frac{1}{r}\exp(-b(r-1)^\beta) \le e^b \frac{1}{r}\exp(-br^\beta),
\end{align*}
where we used the assumptions $r>1$ and $0< \beta \le 1$ for the last inequality. We have proved \eqref{e:ned2}.

Let $\varphi$ be the function \eqref{d:varphi} with $\nu=J_1$ and the dimension $d+2$ (instead of $d$). 
Note that \eqref{e:ned2} implies that $\varphi$  satisfies
$$ c_1 \Phi(r) \le 2c^{-1}\Phi(r/2) =  \frac{c^{-1} r^2}{\int_0^r \frac{s}{\phi(s/2)}ds}\le \varphi(r)=\frac{r^2}{\int_0^r s^{d+3}J_1(s)ds} \le \frac{cr^2}{\int_0^r \frac{s}{\phi(s)}ds} = 2c\Phi(r), r<1. $$
Thus, $J_1$ satisfies \eqref{c:nu1} with $\ell=0$ and \eqref{e:varphi1} since $\Phi$ satisfies \eqref{e:wsc-Phi}. Combining $\varphi(r) \asymp \Phi(r)$ and Proposition \ref{p:hke} for the process $Y$, we have that there is a constant $c_2>0$ satisfying
$$ \frac{-1}{2\pi r}\frac{\partial}{\partial r}p(t,r)=q_t(r) \le c_2 t \G^{(d+2)}(t,r)  \quad
\text{and} \quad \left|\frac{d}{dr} q_t(r)\right| \le c_2t \pinv(t)^{-1} \G^{(d+2)}(t,r) $$
for any $0<t \le T$ and $r>0$. 
From now on, assume $t \in (0,T]$ and $x \in \R^d$. Also, let $r=|x|$. Combining above inequalities and \eqref{e:p-d+2} we have
\begin{align}
\begin{split}
\label{3}
| \nabla^2_x &p(t,x)|=| \frac{\partial^2}{\partial r^2} p(t,r) + \frac{d-1}{r} \frac{\partial}{\partial r} p(t,r)|  =2\pi| \frac{d}{dr} \left( - r q_t(r) + (d-1)q_t(r) \right)| \\ &\le 2\pi d \left(q_t(r) +  r |\frac{d}{dr} q_t(r)  | \right) \le c_3t  \left(1 + r \pinv(t)^{-1} \right)\G^{(d+2)}(t,r) \\
&\le c_4 t \left(1 + r \pinv(t)^{-1} \right)\G_T^{(d+2)}(t,r)
\end{split}
\end{align}
where we used \eqref{e:rho_T} for the last line. Thus, using \eqref{e:rho} we obtain
\begin{align}\label{e:fest1}
 |\nabla_x^2 p(t,x)| \le 2 c_4t  \G^{(d+2)}(t,r) \le 2c_4  \pinv(t)^{-d-2}= 2c_4 \pinv(t)^{-2} \G_T(t,r), \quad  r \le \pinv(t).
 \end{align}
Also, for $\pinv(t)<r \le \pinv(T)$ we have
\begin{align}\begin{split}\label{e:fest2} |\nabla_x^2 p(t,x)| &\le 
 \frac{2c_4t   r} {\pinv(t)} \G^{(d+2)}(t,r) \le \frac{2c_4t   r^2} {\pinv(t)^2} \G^{(d+2)}(t,r) \\ &\le 2c_4\pinv(t)^{-2} \frac{ t}{r^d\Phi(r)}= 2c_4 \pinv(t)^{-2}\G_T(t,r). \end{split}\end{align}
Note that above estimates are valid for any $0<\beta \le 1$. 

Now assume $0<\beta<1$. Let us recall that $J_1$ satisfies \eqref{c:nu1} with $\ell=1$ and \eqref{e:varphi1}. Applying Proposition \ref{p:hke}(c) for the process $Y$ we have
$$ q_t(r) \le c_5 tr^{-1} \exp(-br^\beta) \quad \mbox{and} \quad |\frac{d}{dr} q_t(r)| \le c_5t \pinv(t)^{-1} r^{-1} \exp(-br^\beta) $$
for $r>\pinv(T)$. Thus, by \eqref{3}
\begin{align*}
|\nabla_x^2 p(t,x)| &\le 2\pi d \left(q_t(r) +  r |\frac{d}{dr} q_t(r)  | \right) \\ &\le c_6 t \big(r^{-1}  +  t\pinv(t)^{-1}  \big)\exp(-br^\beta) \le c_7 t\pinv(t)^{-2} \G_T(t,r)
\end{align*}
for $r>\pinv(T)$. Combining this with \eqref{e:fest1}, \eqref{e:fest2} and \eqref{e:rho_T} we obtain
 $$ |\nabla^2_x p_t(x) | \le c_8t \pinv(t)^{-2} \G_T(t,x) \le c_9 \G_{-2}^0(t,x), \quad\quad  0<t \le T, x \in \R^d.$$
 This concludes \eqref{e:grad2} for $0<\beta<1$. 

Similarly, for $\beta=1$ using Proposition \ref{p:hke}(b) for the process $Y$ we have
$$ q_t(r) \le c_{10}t \exp(-\frac{b}{4}r) \quad \mbox{and} \quad  |\frac{d}{dr} q_t(r)| \le c_{10}t \pinv(t)^{-1}  \exp(- \frac{b}{4}r) $$
for $r>\pinv(T)$. Thus,
\begin{align*}|\nabla_x^2 p(t,x)| &\le 2\pi d \left( q_t(|x|) + |x|| \frac{d}{dr} q_t(|x|)| \right)
\le c_{11}t \left( \pinv(t)^{-1}+ r \right) \exp(-\frac{b}{4}r) \\
 &\le c_{12} t\pinv(t)^{-2} \exp(-\frac{b}{5}|x|) \le c_{13} t\pinv(t)^{-2}\G_T(t,r), \quad \quad r>\pinv(T).
 \end{align*}
Hence, combining this with \eqref{e:fest1}, \eqref{e:fest2} and \eqref{e:rho_T} we obtain
$$ |\nabla^2_x p_t(x) | \le c_{14}t \pinv(t)^{-2} \G_T(t,x) \le c_{15}t \G_{-2}^0(t,x), \quad 0<t \le T, x \in \R^d,$$
which is our desired result for $\beta=1$. \qed

\section{Further properties of heat kernel for isotropic L\'evy process}
In this section we assume that $J$ satisfies \eqref{e:J1} with $0<\beta \le 1$ and \eqref{e:J2}, and that nondecreasing function $\phi$ satisfies \eqref{e:lsc-phi} and \eqref{e:phi}. As in the previous section, let $X$ be an istropic unimodal L\'evy process with jumping kernel $J(|y|)dy$ and $p(t,x)$ be the transition density of $X$. Also, let $\sL$ be an infinitesimal generator of $X$.

Recall that  $\delta_f$ is defined in \eqref{d:del}.
The following results are counterpart of \cite[Proposition 3.3]{KSV}.

\begin{prop}\label{t:p}
 For every $T\ge 1$, there exists a constant $0<c=c(d,T,a, a_1,\alpha_1,b,\beta,C_0)$ such that for every $t\in (0,T]$ and $x,y,z\in \R^d$,
	\begin{equation}\label{e:p(a)}
	\left|p(t,x)-p(t,y)\right|\le c\left(\frac{|x-y|}{\pinv(t)}\wedge 1\right)t\left(\G(t,x)+\G(t,y)\right)\, ,
	\end{equation}
	and
	\begin{equation}\label{e:delta(a)}
	\left|\delta_{p}(t,x;z)\right| \le c\left(
	\frac{|z|}{\pinv(t)}\wedge 1\right)^2t\left(\G(t,x\pm z)+\G(t,x)\right),
	\end{equation}
\end{prop}
\pf (a) Since \eqref{e:p(a)} is clearly true when $\pinv(t) \le |x-y|$ by \eqref{e:grad2},  we  assume that $\pinv(t) \ge |x-y|$. Let $\alpha(\theta)=x+\theta(y-x), \,\, \theta \in [0,1]$ be a segment from $x$ to $y$. Then, for any $\theta \in [0,1]$ we have
$$ |\alpha(\theta)| \ge |x| - |x-\alpha(\theta)| \ge |x| - |x-y| \ge |x| - \pinv(t), $$
here we used $|x-y| \le \pinv(t)$ for the last inequality. Thus, we obtain
\begin{align}
|p(t,x)-p(t,y)|&=|\int_0^1 \alpha'(\theta) \cdot \nabla_x p(t, \alpha(\theta))\, d\theta| \le \int_0^1 |\alpha'(\theta)| |\nabla_x p(t, \alpha(\theta))|\, d\theta \nonumber\\
&\le c_1\int_0^1 |\alpha'(\theta)| \pinv(t)^{-1}t\G(t, \alpha(\theta))\, d\theta \le \, c_1|x-y| \pinv(t)^{-1}t\G(t, |x|-\pinv(t)) \nonumber\\
&\le c_2 \, |x-y|\pinv(t)^{-1}t\G(t, x). \nn
\end{align}
Here we used \eqref{e:grad2} with $k=1$ for the second line and \eqref{e:rho-scaling(b)} for the last line. This concludes \eqref{e:p(a)}.

Note that using \eqref{e:grad2} for $k=2$ and following the same argument as the above we can estimate $|\nabla p(t,x)- \nabla p(t,y)|$. Hence, we have a constant $c_3>0$ satisfying
\begin{equation}\label{e2}
|\nabla p(t,x)- \nabla p(t,y) | \le c_3 |x-y| \pinv(t)^{-2} t(\G(t,x)+\G(t,y))
\end{equation}
for $0< t \le T$ and $|x-y| \le \pinv(t)$. 

\noindent
(b)  \eqref{e:delta(a)} is clearly true when $\pinv(t) \le 2|z|$. Now assume $\pinv(t) \ge 2|z|$. Let $\alpha(\theta)=x+\theta z, \,\, \theta \in [-1,1]$ be a segment from $x-z$ to $x+z$. Then, for any $\theta \in [-1,1]$ we have $|\alpha(\theta)| \ge |x| - \pinv(t)/2$, hence
\begin{align}
|\delta_{p}(t,x;z)|
&=|(p(t,x)-p(t,x-z))-(p(t,x+z)-p(t,x))| \nonumber\\
&=\big|\int_0^{1} \alpha'(\theta) \cdot \nabla p(t,\alpha(\theta)) - \alpha'(-\theta) \cdot \nabla p(t,\alpha(-\theta)) d\theta \big|\nonumber \\
&=\big|\int_0^{1} z \cdot (\nabla p(t,\alpha(\theta)) - \nabla p(t,\alpha(-\theta))) d\theta \big|\nonumber \\
&\le 4c_3 \, |z| \pinv(t)^{-2} \left(|z| t \G(t,|x| - \pinv(t))\right) \nn \\
&\le c_4 \, \pinv(t)^{-2} |z|^2 t \G(t,x). \nn
\end{align}
Here we used $|\alpha(\theta)-\alpha(-\theta)| \le 2|z| \le \pinv(t)$ and \eqref{e2} for the first inequality, and \eqref{e:rho-scaling(b)} for the second one.  \qed
\begin{prop}\label{p:int-delta} For every $T\ge 1$, there exist  constants $
	c_i=c_i(d,T, a,a_1,\alpha_1,b,\beta,C_0)
	>0$, $i=1,2$, such that for any $t\in(0,T]$ and $x \in\R^d$,
	\begin{equation}\label{e:int-delta}
	\int_{\R^d} \left|\delta_{p}(t,x;z)\right| J(|z|)\, dz \le c_1 \int_{\R^d}\big( \frac{|z|}{\pinv(t)}\wedge 1\big)^2 t\left(\G(t,x\pm z)+\G(t,x)\right)J(|z|)\, dz \le c_2 \G(t,x)
	\end{equation}
\end{prop}
\pf 
By \eqref{e:delta(a)} we have
\begin{align}
\begin{split}\label{e:3.4.1}
\int_{\R^d} &\left|\delta_{p}(t,x;z)\right| J(|z|)\, dz  \\
&\le c_1 \int_{\R^d}\big( \frac{|z|}{\pinv(t)}\wedge 1\big)^2 t\left(\G(t,x\pm z)+\G(t,x)\right)J(|z|)\, dz \\
&\le c_2 \left(\int_{\R^d}\big( \frac{|z|}{\pinv(t)}\wedge 1\big)^2t\G(t,x + z)J(|z|)\, dz + t\G(t,x)\sP (\pinv(t))
\right) \\
&=: c_2\, (I_1+I_2)\, 
\end{split}
\end{align}
Clearly, by \eqref{e:pp} we have
\begin{equation}\label{e:3.4.2}
I_2 \le c_3 \G(t,x).
\end{equation}
To estimate $I_1$, we divide $I_1$ into two parts as 
\begin{align*}
I_1&=\int_{|z|\le \pinv(t)}\big( \frac{|z|}{\pinv(t)}\big) ^2t \G(t,x+ z)J(|z|)\, dz + \int_{|z|>  \pinv(t)}t \G(t,x+ z)J(|z|)\, dz \\
&=:I_{11}+I_{12}\, . 
\end{align*}
By using \eqref{e:rho-scaling(b)} in the first inequality below and \eqref{e:pp} in the third, we have
\begin{align*}
I_{11}&\le  c_4 t \G(t,x) \int_{|z|\le \pinv(t)}
\big( \frac{|z|^2}{\pinv(t)^2} \land 1\big)J(|z|)\, dz  \\
& \le      c_4 t \G(t,x)  \sP (\pinv(t)) \le c_5 \G(t,x)\, . 
\end{align*}
For the estimates of $I_{12}$, we will use 
\begin{equation}\label{e:J}
 J(|z|) \le c_6 \theta(|z|) = c_6 \G_T(t,z), \quad |z| > \pinv(t),
\end{equation}
which follows from \eqref{e:J1} and \eqref{e:Phi}. Using \eqref{e:rho_T}, \eqref{e:J} and \eqref{e:con}, we arrive
\begin{align*}
I_{12} \le  c_6 a t \int_{|z|>\pinv(t)} \G(t,x-z)\G(t,z)\, dz \le  c_7  \G(t,x).
\end{align*}
Here we used \eqref{e:rho-scaling(a)} for the last inequality. The lemma follows from the estimates of $I_{11}, I_{12}$ and $I_2$. \qed

\subsection{Dependency of $p^\kk$ in terms of $\kk$ }

Recall that 
$$\sL^\kappa f(x) = \lim_{\epsilon \downarrow 0} \int_{ |z| > \epsilon} \left( f(x+z) -f(x) \right) \kappa(x,z) J(|z|)dz$$
where $J: \R_+ \rightarrow \R_+$ is a non-increasing function satisfying \eqref{e:J1} and \eqref{e:J2} with strictly increasing function $\phi$ satisfying \eqref{e:lsc-phi} and \eqref{e:phi}. \\

Let $\kk : \R^d \rightarrow (0,\infty)$ be a symmetric function satisfying
\begin{equation}\label{d:kk}
\kappa_0 \le \kk(z) \le \kappa_1 \quad \mbox{for all } z \in \R^d
\end{equation}
where $\kappa_0$ and $\kappa_1$ are constants in \eqref{c:kap1}.  We denote $Z^\kk$ symmetric L\'evy process  whose jumping kernel is given by $\kk(z)J(|z|), z\in \R^d$. Then the infinitesimal generator of $Z^\kk$ is a self-adjoint operator in $L^2(\R^d)$ and is of the following form:
\begin{align}
\begin{split}\label{d:L^k}
\sL^\kk f(x)&= \lim_{\epsilon \downarrow 0} \int_{ |z| > \epsilon} (f(x+z)-f(x)) \kk(z) J(|z|)dz \\
&= \frac{1}{2} \lim_{\epsilon \downarrow 0} \int_{ |z| > \epsilon} (f(x+z)+f(x-z)-2f(x)) \kk(z)J(|z|)dz. 
\end{split}
\end{align} 
\eqref{e:J1} implies that when $f \in C^2_b(\R^d)$, it is not necessary to take the principal value in the last line in \eqref{d:L^k}. The transition density of $Z^\kk$(i.e., the heat kerenl of $\sL^\kk$) will be denoted by $p^\kk(t,x).$ In this section, we will observe further properties of $p^\kk(t,x).$ \\
	
\begin{remark}\label{r:KSV1} \rm
 The operator \eqref{d:sL} satisfies all conditions in \cite{KSV} with respect to the function $\tG(t,x)$ and $\Phi(r^{-1})^{-1}$ except \cite[(1.7)]{KSV}:
Recall from  Remark \ref{r:KSV} that $\tG(t,x)$ is comparable to the function $\rho(t,x)$ in \cite{KSV}.
Moreover, by Lemma \ref{l:usc-psi}, The characteristic exponent of any symmetric L\'evy process  whose jumping kernel comparable to $J(|z|)$, is comparable to $\Phi(r^{-1})^{-1}$. Clearly \eqref{e:wsc-Phi} and \cite[Remark 5.2]{KSV} with \eqref{e:J1} imply \cite[(1.4), (1.5) and (1.9)]{KSV}.   Also, we obtain gradient estimates with respect to $\tG(t,x)$ in Proposition \ref{p:reg}, which are same as the gradient estimates in \cite[Proposition 3.2]{KSV}. 
Under these observations, one can follow the proofs of \cite{KSV} using \eqref{e:J1} instead of the condition \cite[(1.7)]{KSV} and see that 
\cite[Theorems 1.1--1.3]{KSV} hold under our setting. 
\end{remark}

Using the above Remark \ref{r:KSV1}, from the remainder this paper {\it we use \cite[Theorems 1.1--1.3]{KSV} without any further remark. } \\

Let $\wh\kk:= \kk -\frac{\kappa_0}{2}.$ Then, $\frac{\kappa_0}{2} \le \wh\kk(z) \le \kappa_1$. Let $p^{\wh\kk}$ be the heat kernel of symmetric L\'evy process $Z^{\wh\kk}$ whose jumping kernel is $\wh\kk(z)J(|z|)dz$ and $p^{\frac{\kappa_0}{2}}(t,x)= p(\frac{\kappa_0}{2}t,x)$ be the heat kernel of symmetric L\'evy process $Z^{\frac{\kappa_0}{2}}$ whose jumping kernel is $\frac{\kappa_0}{2}J(|z|)dz$. Without loss of generality, we can assume that $Z^{\wh\kk}$ and $Z^{\frac{\kappa_0}{2}}$ are independent. By \cite[Theorem 1.2]{GKK}, there exists a constant $c=c(T)=c(d,T,a,a_1,\alpha_1,b,\beta,C_0,\kappa_0,\kappa_1)>0$ such that 
\begin{equation}\label{e:p^k}
p^{\wh\kk} (t,x) \le ct \G(t,x) \quad \mbox{for all} \quad 0<t \le T,\,\, x \in \R^d
\end{equation}
for every $\kk$ satisfying \eqref{d:kk}. Also, by Remark \ref{r:KSV1} we have \cite[(3.21)]{KSV}. We record this  for the readers:
\begin{equation}\label{e:p^kk}
\frac{\partial p^\kk(t,x)}{\partial t} =  \sL^\kk p^\kk(t,x), \quad \lim_{t \downarrow 0} p^\kk(t,x) = \delta_0(x).
\end{equation}
\noindent Since $Z^{\wh\kk}$ and $Z^{\frac{\kappa_0}{2}}$ are independent, $Z^\kk$ and $Z^{\wh\kk}+Z^{\frac{\kappa_0}{2}}$ have same distributions. Thus, we have
\begin{align}
\begin{split}\label{e:con-p^kk}
p^\kk(t,x)\quad &= \quad \int_{\R^d} p^{\frac{\kappa_0}{2}}(t,x-y) p^{\wh\kk}(t,y)dy \\
\quad &= \quad \int_{\R^d} p(\frac{\kappa_0}{2}t,x-y) p^{\wh\kk}(t,y)dy.
\end{split}
\end{align}

First we extend  Propositions \ref{p:grad1} and \ref{t:p}--\ref{p:int-delta}. 
\begin{prop}\label{p:p^k}
	There exists a constant $c=c(d,T, a,a_1,\alpha_1,b,\beta,C_0,\kappa_0,\kappa_1)>0$ such that for any $t\in(0,T]$ and $x,y,z\in \R^d, \,$  
	\begin{equation}\label{e:p^k(a)}
	|\nabla_x p^\kk(t,x)| \le ct  \pinv(t)^{-1}\G(t,x)
	\end{equation}
	\begin{equation}\label{e:p^k(b)}
	|p^\kk(t,x)-p^\kk(t,y)| \le ct (\pinv(t)^{-1}|x-y| \land 1 )(\G(t,x)+\G(t,y)) 
	\end{equation}
	\begin{equation}\label{e:p^k(c)}
	|\delta_{p^\kk} (t,x;z)| \le ct ((\pinv(t)^{-1}|z|)^2 \land 1)(\G(t,x\pm z)+\G(t,x)) 
	\end{equation}
	\begin{equation}\label{e:p^k(e)}
	\int_{\R^d} |\delta_{p^\kk}(t,x;z)| J(|z|)dz \le c \G(t,x)
	\end{equation}
\end{prop}
\pf (a) Using \eqref{e:con-p^kk}, \eqref{e:p^kk}, \eqref{e:grad2}, \eqref{e:con} and \eqref{e:rho-scaling(a)} for each line, we obtain
\begin{align*}
|\nabla_x p^\kk (t,x)| &= \Big|\, \nabla_x \int_{\R^d} p(\frac{\kappa_0}{2}t,x-y) p^{\wh\kk}(t,y)dy\, \Big| \\
&\le \Big|\,  \int_{\R^d} \nabla_x p(\frac{\kappa_0}{2}t,x-y) t \G(t,y)dy\,\Big| \\
&\le c_1\int_{\R^d} t \pinv(t)^{-1}\G(\frac{\kappa_0}{2}t,x-y) \times t\G(t,y) dy \\
&\le c_2t\pinv(t)^{-1} \G((1+\frac{\kappa_0}{2})t, x) \\
&\le c_3t\pinv(t)^{-1}\G(t,x).
\end{align*}
	(b) Using \eqref{e:con-p^kk}, \eqref{e:p^kk}, \eqref{e:p(a)}, \eqref{e:con} and \eqref{e:rho-scaling(a)} we obtain
	\begin{align*}
	|p^\kk(t,x)-p^\kk(t,y)| &\le  \int_{\R^d} |p(\frac{\kappa_0}{2},x-z)-p(\frac{\kappa_0}{2},y-z)| p^{\wh\kk}(t,z)dz \\ &\le c_1\int_{\R^d} t (\frac{|x-y|}{\pinv(t)} \land 1)(\G(\frac{\kappa_0}{2}t,x-z)+\G(\frac{\kappa_0}{2}t,y-z))  \times t\G(t,z) dz \\
	&\le c_2t(\frac{|x-y|}{\pinv(t)} \land 1)(\G((1+\frac{\kappa_0}{2})t, x)+\G((1+\frac{\kappa_0}{2})t, y))  \\
	&\le c_3t(\frac{|x-y|}{\pinv(t)} \land 1)(\G(t,x)+\G(t,y)).
	\end{align*}
	(c) We use \eqref{e:con-p^kk}, \eqref{e:p^kk}, \eqref{e:delta(a)}, \eqref{e:con} and  \eqref{e:rho-scaling(a)} for each line to estimate $|\delta_{p^\kk} (t,x;z)|$.
	\begin{align*}
	|\delta_{p^\kk} (t,x;z)| &\le  \int_{\R^d} |\delta_{p} (\frac{\kappa_0}{2}t,x-y;z)| p^{\wh\kk}(t,y)dy \\ &\le c_1\int_{\R^d} t ((\pinv(t)^{-1}|z|)^2 \land 1)(\G(\frac{\kappa_0}{2}t,x-y\pm z)+\G(\frac{\kappa_0}{2}t,x-y))  \times t\G(t,y) dy \\
	&\le c_2t ((\pinv(t)^{-1}|z|)^2 \land 1)(\G((1+\frac{\kappa_0}{2})t, x \pm z)+\G((1+\frac{\kappa_0}{2})t, x))  \\
	&\le c_3t ((\pinv(t)^{-1}|z|)^2 \land 1)(\G(t,x \pm z)+\G(t,x)).
	\end{align*}
	(d) We use \eqref{e:con-p^kk}, Fubini's theorem, \eqref{e:p^kk}, \eqref{e:int-delta}, \eqref{e:con} and  \eqref{e:rho-scaling(a)} for each line to estimate $\int_{\R^d} |\delta_{p^\kk} (t,x;z)|J(|z|)dz$.
	\begin{align*}
	\int_{\R^d} |\delta_{p^\kk} (t,x;z)|J(|z|)dz &\le  \int_{\R^d}\int_{\R^d} |\delta_{p} (\frac{\kappa_0}{2}t,x-y;z)| p^{\wh\kk}(t,y)dy J(|z|)dz \\ &=\int_{\R^d}\left(\int_{\R^d} |\delta_p (\frac{\kappa_0}{2}t,x-y;z)| J(|z|)dz \right)p^{\wh\kk}(t,y)dy \\
	&\le c_1\int_{\R^d}  \G(\frac{\kappa_0}{2}t,x-y)  \times t\G(t,y) dy \\
	&\le c_2 \G((1+\frac{\kappa_0}{2})t, x)  \le c_3 \G(t,x).
	\end{align*}
\qed

Next, we obtain continuity of transition density with respect to $\kk$. This is the counterpart of \cite[Theorem 3.5]{KSV}.

\begin{thm}	There exists a constant $c=c(d,T,a, a_1,\alpha_1,b,\beta,C_0,\kappa_0,\kappa_1)>0$ such that
	for any  two symmetric functions $\kk_1$ and $\kk_2$ in $\R^d$ satisfying \eqref{d:kk}, any $t\in(0,T]$ and $x\in \R^d$, we have 
	\begin{equation}
	\left|p^{\kk_1}(t,x)-p^{\kk_2}(t,x)\right| \le c \|\kk_1-\kk_2\|_{\infty}\ t\G(t,x)\, , \label{e:p^k12(a)}
	\end{equation}
	\begin{equation}
	\left|\nabla p^{\kk_1}(t,x) -\nabla p^{\kk_2}(t,x)\right| \le   c  \|\kk_1-\kk_2\|_{\infty}\pinv(t)^{-1} t\G(t,x)\label{e:p^k12(b)} 
	\end{equation}
	and
	\begin{equation}\label{e:p^k12(c)}
	\int_{\R^d} | \delta_{p^{\kk_1}} (t,x;z) - \delta_{p^{\kk_2}}(t,x;z)| J(|z|)dz \le c \|\kk_1-\kk_2\|_{\infty}\G(t,x).
	\end{equation}
\end{thm}

\pf (a) $p^{\kk_1}(s,y)$ is uniformly bounded on $s \in [t/2,t]$ by \eqref{e:p^k} and $\displaystyle \lim_{s \rightarrow t} p^{\kk_2}(t-s,x-y)= \delta_0(x-y)$ by \eqref{e:p^kk}. Thus, we have $$\displaystyle \lim_{s \uparrow t} \int_{\R^d}p^{\kk_1}(s,y)p^{\kk_2}(t-s,x-y)\, dy= p^{\kk_1}(t,x).$$  By the similar way, we get $$\displaystyle \lim_{s \downarrow 0} \int_{\R^d}p^{\kk_1}(s,y)p^{\kk_2}(t-s,x-y)\, dy= p^{\kk_2}(t,x).$$
Hence, for $t \in (0,T]$ and $x \in \R^d$,
$$ |p^{\kk_1}(t,x) - p^{\kk_2}(t,x)| = \left|\int_0^t \frac{d}{ds}\left(\int_{\R^d}p^{\kk_1}(s,y)p^{\kk_2}(t-s,x-y)\, dy\right) ds \right|. $$
Using \eqref{e:p^kk} in the second line, the fact that $\LL^{\kk_1}$ is self-adjoint in the third line and \eqref{d:L^k} in the fourth line, we have
\begin{align*}
 \lefteqn{\int_0^{t/2} \frac{d}{ds}\left(\int_{\R^d}p^{\kk_1}(s,y)p^{\kk_2}(t-s,x-y)\, dy\right) ds}\\
	&=\int_0^{t/2} \left(\int_{\R^d}\left(\LL^{\kk_1}p^{\kk_1}(s,y) p^{\kk_2}(t-s,x-y)-p^{\kk_1}(s,y)\LL^{\kk_2}p^{{\kk_2}}(t-s,x-y)\right)dy\right)ds\\
	&=\int_0^{t/2} \left(\int_{\R^d} p^{\kk_1}(s,y)\left(\LL^{\kk_1}-\LL^{\kk_2}\right) p^{\kk_2}(t-s,x-y)dy\right)ds\\
	&=\frac12 \int_0^{t/2} \left(\int_{\R^d} p^{\kk_1}(s,y)\left(\int_{\R^d}\delta_{p^{\kk_2}}(t-s,x-y; z)(\kk_1(z)-\kk_2(z))J(|z|)dz\right)dy\right)ds.
\end{align*}
Hence, by using \eqref{e:p^k(e)},
\eqref{e:p^k(a)} 
and the convolution inequality \eqref{e:con}, we have
\begin{align*}
\lefteqn{ \left|\int_0^{t/2} \frac{d}{ds}\left(\int_{\R^d}p^{\kk_1}(s,y)p^{\kk_2}(t-s,x-y)\, dy\right) ds \right|}\\
	& \le \frac12 \|\kk_1-\kk_2\|_{\infty} \int_0^{t/2}\left(\int_{\R^d}p^{\kk_1}(s,y)\left(\int_{\R^d}\left|\delta_{p^{\kk_2}}(t-s,x-y; z)\right|J(|z|)dz\right)dy\right)ds\\
	&\le c_1 \|\kk_1-\kk_2\|_{\infty} \int_0^{t/2} \int_{\R^d} s \G(s,y)\G(t-s,x-y)\, dy\, ds\\
	&\le c_2\|\kk_1-\kk_2\| _{\infty} \int_0^{t/2} s \big( s^{-1} + (t-s)^{-1} \big) \G(t,x) ds \le c_3 \|\kk_1-\kk_2\|t\G(t,x), 
\end{align*}
for all $t \in (0, T]$ and $x \in \R^d$. By the similar way, we also obtain 
\begin{align*}
\lefteqn{\left|\int_{t/2}^{t} \frac{d}{ds}\left(\int_{\R^d}p^{\kk_1}(s,y)p^{\kk_2}(t-s,x-y)\, dy\right) ds \right|}\\
&=\frac12 \left|\int_{t/2}^{t} \left(\int_{\R^d} p^{\kk_2}(t-s,y)\left(\int_{\R^d}\delta_{p^{\kk_1}}(s,x-y; z)(\kk_1(z)-\kk_2(z))J(|z|)dz\right)dy\right)ds  \right|\\
&\le c_1 \|\kk_1-\kk_2\|_{\infty} \int_{t/2}^{t} \int_{\R^d} (t-s) \G(s,y)\G(t-s,x-y)\, dy\, ds \\
&\le c_3\|\kk_1-\kk_2\|_{\infty}t \G(t,x).
\end{align*}
Therefore, we arrive
\begin{align*}
\lefteqn{|p^{\kk_1}(t,x)-p^{\kk_2}(t,x) | } \\
& \le \left| \int_0^{t/2} \frac{d}{ds}\left(\int_{\R^d}p^{\kk_1}(s,y)p^{\kk_2}(t-s,x-y)\, dy\right) ds \right| + \left|\int_{t/2}^t \frac{d}{ds}\left(\int_{\R^d}p^{\kk_1}(s,y)p^{\kk_2}(t-s,x-y)\, dy\right) ds  \right|\\
&\le 2c_3\|\kk_1-\kk_2\|_{\infty} t\G(t,x). 
\end{align*}
\noindent 
(b) Set $\wh\kk_i(z):=\kk_i(z)-{\kappa_0}/{2}, \,\, i=1,2.
$
Using \eqref{e:con-p^kk},  \eqref{e:grad2},  \eqref{e:p^k12(a)}, \eqref{e:con} and 
\eqref{e:rho-scaling(a)}, we have  that for all  $t \in (0, T]$ and   $x \in \R^d$,
\begin{align*}
\left|\nabla p^{\kk_1}(t,x) -\nabla p^{\kk_2}(t,x)\right|  &=  \left|\int_{\R^d}\nabla p\left(\frac{\kappa_0}{2}t, x-y\right) \left(p^{\wh\kk_1}(t,y)-p^{\wh\kk_2}(t,y)\right)dy\right|\\
	& \le   c_1\|\kk_1-\kk_2\|_{\infty} \pinv(t)^{-1}t^2 \int_{\R^d}\G(\frac{\kappa_0}{2}t,x-y)\G(t,y)\, dy\\
	&\le  c_2 \|\kk_1-\kk_2\|_{\infty} \pinv(t)^{-1}t^2  t^{-1}\G((1+\frac{\kappa_0}{2})t,x)\\
	&\le  c_3 \|\kk_1-\kk_2\|_{\infty} \pinv(t)^{-1} t\G(t,x)\, .
\end{align*}
(c) By using \eqref{e:con-p^kk}, \eqref{e:delta(a)}, \eqref{e:p^k12(a)}, \eqref{e:con} and \eqref{e:rho-scaling(a)} we  have that for any $t \in (0,T]$ and $x,z \in \R^d$,
\begin{align*}
\lefteqn{|\delta_{p_{\kk_1}}(t,x;z) - \delta_{p_{\kk_2}}(t,x;z) | = \left| \int_{\R^d} \delta_p(\frac{\kappa_0}{2}t, x-y;z) \left(p^{\wh\kk_1}(t,y) - p^{\wh\kk_2}(t,y) \right) dy \right| }\\
		& \le  c_1 \|\kk_1-\kk_2\|_{\infty} \big( \pinv(t)^{-1}|z| \land 1 \big)^2 t^2 \int_{\R^d} (\G(\frac{\kappa_0}{2}t,x-y \pm z) + \G(t,x-y) \G(t,y)dy \\
		&\le c_2 \|\kk_1-\kk_2\|_{\infty} \big( \pinv(t)^{-1}|z| \land 1 \big)^2 t^2 \big( t^{-1} \G((1+\frac{\kappa_0}{2})t,x \pm z) + \G((1+\frac{\kappa_0}{2})t,x) \big) \\
		&\le c_3 \|\kk_1-\kk_2\|_{\infty}  \big( \pinv(t)^{-1}|z| \land 1 \big)^2 t \big( \G(t,x \pm z) + \G(t,x) \big).
		\end{align*}
Integrating above inequality we obtain that 
\begin{align*}
\lefteqn{\int_{\R^d} | \delta_{p^{\kk_1}} (t,x;z) - \delta_{p^{\kk_2}}(t,x;z)| J(|z|)dz} \\ & \le c_3t \|\kk_1-\kk_2\|_{\infty} \int_{\R^d} \big( \frac{|z|}{\pinv(t)} \land 1 \big)^2 \big( \G(t,x \pm z) + \G(t,x) \big) J(|z|)dz \\ &
\le c_4 \|\kk_1-\kk_2\|_{\infty} \G(t,x),
\end{align*}
where the last inequality follows from Proposition \ref{p:int-delta}. \qed

Estimates in this section are almost same with \cite[Section 2 and 3]{KSV} except these: First of all, the function $\G$ is different from \cite{KSV}, hence our estimates are more precise  than estimates in \cite{KSV}. However, 
 we don't have estimates for third derivatives in terms of $\G$ of the heat kernel in Proposition \ref{p:grad1}. Thus, we do not have the improvements on 
 \cite[(3.14) and (3.18)]{KSV}, which are used for the gradient estimate of the function $p^\kappa(t,x,y)$ in Theorems \ref{t:main1}-\ref{t:main3}, for instance, \cite[Theorem 1.1(2) and 1.2(4)]{KSV}.  
\section{Estimates of $p^\kappa(t,x,y)$
}
For the remainder of this paper, we always assume that $\kappa:\R^d\times \R^d\to (0,\infty)$ is a Borel function 
satisfying \eqref{c:kap1} and \eqref{c:kap2}, that $J$ satisfies \eqref{e:J1}-\eqref{e:J2} with the function $\phi$ satisfying \eqref{e:lsc-phi} and \eqref{e:phi}.

 For a fixed $y\in \R^d$, let $\kk_y(z)=\kappa(y,z)$ and let $\sL^{\kk_y}$ be the freezing operator defined by
\begin{align}
\label{d:fre}
\LL^{\kk_y}f(x)=
\lim_{\eps \downarrow 0}\int_{|z| > \eps} \delta_f(x;z)  \kappa(y,z)J(|z|)dz.
\end{align}
Let $p_y(t,x):=p^{\kk_y}(t,x)$ be the heat kernel of the operator $\LL^{\kk_y}$.  
Note that $\kk_y$ satisfies \eqref{d:kk} so that there exists a constant $c>0$ such that
\begin{equation}\label{e:p_y}
p_y(t,x) \le c t \G(t,x) \quad \mbox{for all } x,y\in\R^d, t \in (0,T].
\end{equation}

By Remark \ref{r:KSV1} and \cite[Theorem 1.1]{KSV}, 
we have a continuous function $p^\kappa(t,x,y)$ on $(0,\infty) \times \R^d \times \R^d$ solving \eqref{e:main1-1} and it satisfies
$$ p^\kappa(t,x,y) \le c t \tG(t,x-y), \quad \quad 0<t \le T \, \mbox{ and } \, x \in \R^d $$

In this section, we will investigate further estimates and regularity of   $p^\kappa(t,x,y)$. We first recall the construction of $p^\kappa$ from  \cite[section 4]{KSV}. For $t>0$ and $x,y \in \R^d$, define
\begin{equation}\label{d:q0}
q_0(t,x,y):= \frac 12 \int_{\R^d} \delta_{p_y} (t,x-y;z)(\kappa(x,z)-\kappa(y,z)) J(|z|)dz = \big(\sL^{\kk_x} -\sL^{\kk_y}\big) p_y(t,\cdot) (x-y).
\end{equation}
For $n\in \N$, we inductively define the function $q_n(t,x,y)$ by
\begin{equation}\label{d:qn}
q_n(t,x,y):= \int_0^t \int_{\R^d} q_0(t-s,x,z)q_{n-1}(s,z,y) dzds
\end{equation}
and
\begin{equation}\label{d:q}
q(t,x,y):= \sum_{n=0}^\infty q_n(t,x,y).
\end{equation}
Finally we define 
\begin{equation}\label{d:phi-t}
\phi_y(t,x):= \int_0^t \phi_y(t,x,s)ds = \int_0^t \int_{\R^d} p_z(t-s,x-z)q(s,z,y)dzds
\end{equation}
and 
\begin{equation}\label{d:p^kap}
p^\kappa (t,x,y) := p_y(t,x-y) + \phi_y(t,x) = p_y(t,x-y) + \int_0^t \int_{\R^d} p_z(t-s,x-z)q(s,z,y)dzds.
\end{equation}
As \cite[section 4]{KSV}, the definitions in \eqref{d:q0}--\eqref{d:p^kap} are well-defined. In other words, each integrand in \eqref{d:q0}--\eqref{d:p^kap} is integrable and series in \eqref{d:q} is absolutely converge on $(0,\infty) \times \R^d \times \R^d$.

In the next lemma, we will establish the upper bounds of $p^\kappa$.

\begin{thm}\label{t:p^kap}
	 For every $T \ge 1$ and $\delta_0 \in (0,\delta] \cap (0,\alpha_1/2)$, there are constants $c_1$ and $c_2$ such that for any $t \in (0,T]$ and $x,y \in \R^d$, 
	\begin{equation}\label{e:ph}
	|\phi_y(t,x)| \le c_1 t \big(\G_0^{\delta_0} + \G_{\delta_0}^0\big)(t,x-y)
	\end{equation}
	and
	\begin{equation}\label{e:p^kap}
	p^\kappa(t,x,y) \le c_2 t \G(t,x-y).
	\end{equation}
	The constants $c_1$ and $c_2$ depend on $d,T,a, a_1,\alpha_1,b,\beta,C_0,\delta_0,\delta,\kappa_0,\kappa_1$ and $\kappa_2$.
\end{thm}
\pf 
\noindent We first claim that for $n \in \N_0$, 
\begin{equation} \label{e:qn}
|q_n(t,x,y)| \le d_n \big( \G_{(n+1){\delta_0}}^0 + \G_{n{\delta_0}}^{\delta_0} \big)(t,x-y)
\end{equation}
with $$d_n:=(16C(\delta_0,T) c_2)^{n+1} \prod_{k=1}^n B({\delta_0}/2, k{\delta_0}/2)= (16C c_2)^{n+1} \frac{\Gamma({\delta_0}/2)^{n+1}}{\Gamma((n+1){\delta_0}/2)}$$
where $C=C(\delta_0,T)$ is the constant in \eqref{e:con(b)}. Without loss of generality, we  assume that $C \ge 1/16$. \\

\noindent For $n=0$, using \eqref{d:q0},
\eqref{c:kap1}, \eqref{c:kap2} and \eqref{e:p^k(e)} we have
\begin{align*}
|q_0(t,x,y)| &\le  \frac 12 \int_{\R^d} |\delta_{p_y} (t,x-y;z)(\kappa(x,z)-\kappa(y,z))| J(|z|)dz \nn \\
&\le c_1  \big(|x-y|^{\delta_0} \land 1 \big) \int_{\R^d} |\delta_{p_y} (t,x-y;z)| J(|z|)dz \nn \\
&\le c_2 \big(|x-y|^{\delta_0} \land 1 \big) \G(t,x-y) = c_2 \G_0^{\delta_0} (t,x-y).
\end{align*}

\noindent Suppose that \eqref{e:qn} is valid for $n$. Then for $t \le T$,
\begin{align*}
|q_{n+1}(t,x,y)| &\le \int_0^t \int_{\R^d} |q_0(t-s,x,z)q_n(s,z,y)|dzds \\
&\le c_2 d_n \int_0^t \int_{\R^d} \G_0^{\delta_0} (t-s,x-z) \big( \G^0_{(n+1){\delta_0}} + \G^{\delta_0}_{n{\delta_0}} \big) (x,z-y) dzds \\
&\le 16C c_2 d_n B({\delta_0}/2, (n+1){\delta_0}/2) \big( \G^0_{(n+2){\delta_0}} + \G^{\delta_0}_{(n+1){\delta_0}} \big)(t,x-y) \\&= d_{n+1} \big( \G^0_{(n+2)\delta} + \G^{\delta_0}_{(n+1){\delta_0}} \big)(t,x-y)
\end{align*}
here we used induction in the second line, and used \eqref{e:con(c)} and \eqref{e:con(d)} in the last line. For the third line, we need the following: let $\theta=\eta=1$, $\gamma_1=\delta_2=0$, $\delta_1={\delta_0}$ and $\gamma_2=(n+1){\delta_0}$ which satisfy conditions in Lemma \ref{l:con}(c) since ${\delta_0} \in  (0,\alpha_1/2)$. Then, by \eqref{e:con(c)} we have
\begin{align*}
&\int_0^t \int_{\R^d} \G_0^{\delta_0} (t-s,x-z)  \G^0_{(n+1){\delta_0}}(s,z-y) dzds \\
&\le 4C B({\delta_0}/2, (n+1){\delta_0}/2) \big(\G^0_{(n+2){\delta_0}} + \G^{\delta_0}_{(n+1){\delta_0}}+\G^0_{(n+2){\delta_0}}\big) (t,x-y)  \\
&\le 8C B({\delta_0}/2, (n+1){\delta_0}/2) \big( \G^0_{(n+2){\delta_0}} + \G^{\delta_0}_{(n+1){\delta_0}} \big)(t,x-y).
\end{align*}
Also, letting $\theta=\eta=1$, $\gamma_1=0$, $\delta_1=\delta_2={\delta_0}$ and $\gamma_2={\delta_0}$ which satisfy conditions in Lemma \ref{l:con}(c),
\begin{align*}
&\int_0^t \int_{\R^d} \G_0^{\delta_0} (t-s,x-z)  \G^{\delta_0}_{n{\delta_0}}(x,z-y) dzds \\
&\le 4C B({\delta_0}/2, (n+1){\delta_0}/2) \big(\G^0_{(n+2){\delta_0}} + \G^{\delta_0}_{(n+1){\delta_0}}+\G^{\delta_0}_{(n+1){\delta_0}}\big) (t,x-y)  \\
&\le 8C B({\delta_0}/2, (n+1){\delta_0}/2) \big( \G^0_{(n+2){\delta_0}} + \G^{\delta_0}_{(n+1){\delta_0}} \big)(t,x-y).
\end{align*}

\noindent Thus, \eqref{e:qn} is valid for all $n \in \N_0$. Note that
\begin{equation}
\label{31}
\sum_{n=0}^\infty d_n \pinv(T)^{{\delta_0}} := C_1(\delta_0,T) <\infty
\end{equation}
since $ \frac{d_{n+1} \pinv(T)^{(n+1){\delta_0}}}{d_n \pinv(T)^{n{\delta_0}}} = 16Cc_2 \pinv(T)^{\delta_0} B({\delta_0}/2, (n+1){\delta_0}/2) \rightarrow 0 $ as $n \rightarrow \infty$. So, by using \eqref{e:rho-gam} in the second line we obtain
\begin{align*}
&\sum_{n=0}^\infty |q_n(t,x,y)| \le \sum_{n=0}^\infty d_n \big( \G^0_{(n+1){\delta_0}} + \G^{\delta_0}_{n\delta_0} \big) (t,x-y)  \nn \\
&\le \sum_{n=0}^\infty d_n \pinv(T)^{n\delta_0} \big( \G^0_{\delta_0} + \G_0^{\delta_0}) (t,x-y) = C_1 \big( \G^0_{\delta_0} + \G_0^{\delta_0}) (t,x-y)
\end{align*} 
for $t \le T$. Therefore, for every $t \in (0,T]$ and $x,y \in \R^d$,
\begin{equation}
\label{e:q}
|q(t,x,y)| \le C_1 \big( \G^0_{\delta_0} + \G_0^{\delta_0}) (t,x-y).
\end{equation} 
To obtain \eqref{e:ph} and \eqref{e:p^kap}, we calculate that 
\begin{align}
|\phi_y(t,x)| &\le \int_0^t \int_{\R^d}p_z(t-s,x-z) |q(s,z,y)|\, dz\, ds \nn \\
&\le  c_3\int_0^t \int_{\R^d} (t-s)\G(t-s,x-z)\left(\G_{\delta_0}^0+\G_0^{\delta_0}\right)(s,z-y)\, dz\, ds \nn \\
&\le  c_4 t\left(\G_{\delta_0}^0+\G_0^{\delta_0}\right)(t,x-y) \nn \\ 
&\le  2c_4 \pinv(T)^{\delta_0} t\G(t,x-y)= c_5 t \G(t,x-y), \quad \text{ for all } t\in (0,T]\, .\nn
\end{align}
Here we used \eqref{e:p^kk} and \eqref{e:q} for the second line, \eqref{e:con(c)} for the third line and \eqref{e:rho1} for the last line.
Therefore, using \eqref{e:p^kk} we arrive $p^{\kappa}(t,x,y)\le p_y(t,x-y)+|\phi_y(t,x)|\le c_6 t \G(t,x-y)$. \qed

We concludes this section with some fractional estimates on $p^\kappa(t,x,y)$. 

\begin{lemma}\label{l:kap} For every $T \ge 1$ and $\gamma \in (0,1] \cap (0,\alpha_1)$, there exists a constant $c_3$
	such that for any $t\in (0, T]$ and $x, x', y\in \R^d$,
	$$
	\left|p^{\kappa}(t,x,y)-p^{\kappa}(t,x',y)\right| \le c_3 |x-x'|^{\gamma}\ t\Big(\G_{-\gamma}^0(t,x-y)+\G_{-\gamma}^0(t,x'-y)\Big).
	$$
The constant $c_3$ depends on $d,T,a, a_1,\alpha_1,b,\beta,C_0,\gamma,\delta,\kappa_0,\kappa_1$ and $\kappa_2$.
\end{lemma}
\pf 
Assume that $x, x', y\in \R^d$ and $t\in (0, T]$. By \eqref{e:p^k(b)} and the fact that $\gamma\le 1$, we have
\begin{align}\label{e:5.3}
\begin{split}
	|p_z(s,x-z)-p_z(s,x'-z)|
	&\le c_1 |x-x'|^{\gamma} s \Phi^{-1}(s)^{-\gamma}\big(\G(s,x-z)+\G(s,x'-z)\big)\\
	&\le  c_1 |x-x'|^{\gamma} s \big(\G_{-\gamma}^0(s,x-z)+\G_{-\gamma}^0(s,x'-z)\big)\, .
\end{split}
\end{align}
for any $0<s \le T$ and $z \in \R^d$. Thus, by \eqref{e:q}, \eqref{e:5.3} and a change of the variables, we further have that for $\delta_0 := (\delta \land \alpha_1/4) \in (0,\delta] \cap (0,\alpha_1/2)$,
\begin{align*}
&|\phi_y(t,x)-\phi_y(t,x')|
\le  \int_0^t \int_{\R^d}|p_z(t-s,x-z)-p_z(t-s,x'-z)| \, |q(s,z,y)|\, dz \, ds\\
\le & c_2 |x-x'|^{\gamma}\int_0^t \int_{\R^d}(t-s) 
\Big(\G_{-\gamma}^0(t-s,x-z)+\G_{-\gamma}^0(t-s,x'-z)\Big)
\big(\G_0^{\delta_0}+\G_{\delta_0}^0\big)(s,z-y)\, dz\, ds\\
\le & c_3 |x-x'|^{\gamma} t
\Big(\G_{-\gamma+\delta_0}^0(t,x-y)+\G_{-\gamma}^{\delta_0}(t,x-y)+\G_{-\gamma+\delta_0}^0(t,x'-y)+\G_{-\gamma}^{\delta_0}(t,x'-y)\Big)\\
\le & 2c_3 \Phi^{-1}(T)^{\delta_0} 
|x-x'|^{\gamma} t \big(\G_{-\gamma}^0(t,x-y)+\G_{-\gamma}^0(t,x'-y)\big), \quad \text{ for all } t\in (0,T] \, .
\end{align*}
Since $\gamma < \alpha_1$, the penultimate line follows from \eqref{e:con(c)} (with $\theta=0$), 
and the last line by \eqref{e:rho-gam} and \eqref{e:rho-del}. The lemma follows by combining above two estimates and \eqref{d:p^kap}. \qed
\section{Proof of Theorems \ref{t:main1}-\ref{t:main3}}

In this section we prove the main theorems in Section 1. We first prove that the function $p^\kappa(t,x,y)$ defined by \eqref{d:p^kap} satisfies all statements in Theorems \ref{t:main1}-\ref{t:main3},
then we show that $p^\kappa(t,x,y)$ is the unique solution to \eqref{e:main1-1} satisfying (i)--(iii) in Theorem \ref{t:main1}.

\textbf{Proof of Theorems \ref{t:main2} and \ref{t:main3}.}
It follows from
Remark \ref{r:KSV1} that we can apply the results in \cite[Theorem 1.1-1.4]{KSV} for operator \eqref{d:sL} with the function $\tG(t,x)$. Note that the function $p^\kappa(t,x,y)$ in \cite[Theorems 1.1-1.4]{KSV} is constructed by the same way as \eqref{d:p^kap}. Therefore, Theorems \ref{t:main2} and \ref{t:main3} except \eqref{e:main2-3} immediately follow from Remarks \ref{r:KSV} and  \ref{r:KSV1}, and \cite[Theorem 1.1(iii), 1.2 and 1.3]{KSV}.  Finally \eqref{e:main2-3} is proved in Lemma \ref{l:kap}.\qed

Now we prove the lower bound estimates in Theorem \ref{t:main1} and Corollary \ref{c:tshke} for the function $p^\kappa(t,x,y)$ in \eqref{d:p^kap}. By Theorems \ref{t:main2} and \ref{t:main3}, we have that $(P_t^{\kappa})_{t \ge 0}$ defined by $p^\kappa(t,x,y)$ in \eqref{d:p^kap} with \eqref{e:intro-semigroup} is a Feller semigroup and there exists a Feller process $X=(X_t, \P_x)$ corresponding to $(P_t^{\kappa})_{t  \ge 0}$. Moreover, by \eqref{e:main3-1} for $f\in C_b^{2,\eps}(\R^d)$, 
\begin{align}
\label{e:MG}
f(X_t)-f(x)-\int_0^t \LL^{\kappa}f(X_s)\, ds
\end{align}
is a martingale with respect to the filtration $\sigma(X_s, s \le t)$.
Therefore, by the same argument as that in \cite[Section 4.4]{CZ16}, we have 
the following L\'evy system formula:
for every function $f:\R^d \times \R^d \to [0,\infty)$ vanishing on the diagonal and every stopping time $S$,
\begin{align}
\label{e:LSF}
\E_x \sum_{0<s\le S} f(X_{s-}, X_s)=\E_x \int_0^S f(X_s,y)J_X(X_s,dy) ds\, ,
\end{align}
where $J_X(x,y):=\kappa(x,y-x)J(|x-y|)$. For $A \in \sB(\R^d)$
we define
$\tau_A:=\inf\{t \ge 0: X_t \notin A\}
$ be the exit time from $A$.

Using \eqref{e:MG}, \eqref{c:kap1} and  \eqref{e:pp}, the proof of the following result is the same as the one in \cite[Lemma 5.7]{KSV}.
We skip the proof. 
\begin{lemma}\label{l:ep}
	Let $T \ge 1$. For each $\eps \in (0,1)$ there exists $\lam=\lam(\eps)>0$ such that for  every $0<r\le \pinv(T)$,
	\begin{equation}\label{e:ep}
	\sup_{x \in \R^d} \P_x\left(\tau_{B(x,r)}\le \lam \Phi(r) \right)\le \eps\, .
	\end{equation}
\end{lemma}
We record that by \eqref{e:ep}, for any $x \in \R^d$ and $0<r \le \pinv(T)$ we have
\begin{equation}\label{e:exit}
\E_x[\tau_{B(x,r)}] \ge \lam(1/2) \Phi(r) \P_x(\tau_{B(x,r)} > \lam \Phi(r)) \ge \frac{\lam}{2} \Phi(r) = c\Phi(r).
\end{equation}

Now we are ready to prove the lower bound in \eqref{e:main4-1}. 

\begin{lemma}\label{l:main2} The function $p^\kappa(t,x,y)$ in \eqref{d:p^kap} satisfies \eqref{e:main4-1}.
	\end{lemma}
\pf {Fix $T \ge 1$. Let $p_y(t,x)$ be the heat kernel of the freezing operator in \eqref{d:fre}, and $J_y(z):=\kappa(y,z)J(|z|)$ and $\psi_y(z)$ be the corresponding L\'evy measure and characteristic exponent, respectively. By \cite[Theorem 2]{KS15}, there exist constants $C_1,C_2>0$ such that 
\begin{equation}\label{e:lower}
p_y(t,x)\ge C_1 \pinv(t)^{-d}, \quad \quad t \in (0,T], y \in \R^d \,\, \mbox{and} \,\, |x| \le C_2\pinv(t).
\end{equation}
Indeed, $J_y(z)dz$ is symmetric and infinite L\'evy measure by \eqref{c:kap1} and that $J(|z|)dz$ is infinite L\'evy measure. To check the condition \cite[(3)]{KS15}, we need to show that there exists a constant $c>0$ such that}
\begin{align}
\label{e:nn0}
\int_{\R^d} e^{-t \psi_y(z)} |z|dz \le ch_y(t)^{-d-1}, \quad 0<t, y \in \R^d 
\end{align}
where $h_y(t):= \frac{1}{\Psi_y^{-1}(t^{-1})}$ and $\Psi_y(r):= \sup_{|z| \le r} \psi_y(z)$. Let $\PP(r):=\int_{\R^d} \big(1 \land \frac{|z|^2}{r^2} \big) J(|z|)dz$ and $\PP_y(r):= \int_{\R^d} \big(1 \land \frac{|z|^2}{r^2} \big) J_y(|z|)dz$. Then, by \cite[(11)]{KS15} we have
\begin{align}
\label{e:nn1}
c_1\kappa_0 \PP(r^{-1}) \le c_1 \PP_y(r^{-1}) \le \Psi_y(r) \le 2\PP_y(r^{-1}) \le 2\kappa_1 \PP(r^{-1}), \quad r>0. 
\end{align}
 On the other hand, by the symmetry of $J_y$ and \cite[(10)]{KS15}, we have
\begin{align*}
 \psi_y(z) &\ge (1-\cos 1)\int_{|\xi| \le 1/|z|} |\xi \cdot z|^2 J_y(d\xi) \ge \kappa_0(1-\cos 1)  \int_{|\xi| \le 1/|z|} |\xi \cdot z|^2 J(|\xi|)d\xi. 
\end{align*}
Since by a rotation
$$\int_{|\xi| \le 1/|z|} |\xi \cdot z|^2 J(|\xi|)d\xi= |z|^2 \int_{|\xi| \le 1/|z|} \xi_i^2J(|\xi|)d\xi, \quad i=1, \dots, d,$$
we have
\begin{align*}
\psi_y(z)  \ge d^{-1} \kappa_0(1-\cos 1)   |z|^2 \int_{|\xi| \le 1/|z|} |\xi|^2J(|\xi|)d\xi.
\end{align*}
Thus, when $|z| \le 1$ we have
$$ \psi_y(z) \ge d^{-1} \kappa_0(1-\cos 1)|z|^2 \int_{|\xi| \le 1} |\xi|^2J(d\xi) \ge c_2|z|^2 = c_3\Phi(|z|^{-1}), $$
whereas by \eqref{e:J1} we have
$$ \psi_y(z) \ge d^{-1} \kappa_0(1-\cos 1)|z|^2 \int_{|\xi| \le 1/|z|} |\xi|^2 J(d\xi) \ge c_4 |z|^2 \int_0^{1/|z|} \frac{s}{\phi(s)}ds = c_4\Phi(|z|^{-1}) $$
for $|z| \ge 1$. Therefore, using \eqref{e:pp} and \eqref{e:pru} we obtain
\begin{align}
\label{e:nn2}
\psi_y(z) \ge c_5 \Phi(|z|^{-1}) \ge c_6 \PP(|z|) \ge (c_6/2)\psi(|z|).
\end{align} 
Moreover,  \eqref{e:pp} and \eqref{e:nn1} also imply that $h_y(t) \asymp \pinv(t) \asymp h(t):= \frac{1}{\Psi^{-1}(t^{-1})}$. From this and  and \eqref{e:nn2} we can follow the proof of Lemma \ref{l:E} and obtain \eqref{e:nn0} as 
\begin{align}
\label{e:nn3}
\int_{\R^d} e^{-t  \psi_y(z)} |z|dz \le \int_{\R^d} e^{-c_6t  \psi(z)/2} |z|dz \le 
 c_7h(t)^{-d-1} \le c_8h_y(t)^{-d-1}, \quad 0<t, y \in \R^d.
\end{align}
 Note that every constant above is independent of $y$. Therefore, letting $f(r) \equiv 0$ we obtain all conditions in \cite[Theorem 2]{KS15} so we have \eqref{e:lower} where $C_1>0$ is independent of $y$. 
 
The rest of  the proof is almost identical to the one of \cite[Theorem 1.4]{KSV}. 
Note that there is minor gap in  \cite[(5.36)]{KSV}.
We provide the full details here including the correction of  \cite[(5.36)]{KSV}.

Choose $t_0 \in (0,T]$ small enough to satisfy $2c_1 \pinv(t_0)^{\delta_0} \le C_1 /2$ where $c_1$ and $\delta_0$ are constants in \eqref{e:ph}. Then, using \eqref{e:ph} and \eqref{e:rho1} we have that for any $0<t \le t_0$ and $x,y \in \R^d$,
$$ |\phi_y(t,x)| \le c_1t \big(\G_0^{\delta_0} + \G^0_{\delta_0} \big)(t,x-y) \le 2c_1 \pinv(t_0)^{\delta_0}  t\G(t,x-y) \le 2c_1 \pinv(t_0)^{\delta_0}  \pinv(t)^{-d} \le \frac{C_1}{2} \pinv(t)^{-d}.$$
Thus, combining above inequality and \eqref{e:lower} we obtain
$$p^\kappa(t,x,y) = p_y(t,x-y) + \phi_y(t,x) \ge p_y(t,x-y) - |\phi_y(t,x)| \ge \frac{C_1}{2}\pinv(t)^{-d}$$
for $0<t \le t_0$ and $x,y \in \R^d$ with $|x-y| \le C_2\pinv(t)$. By \eqref{e:main2-2} and iterating at most $n_0:=\lfloor T/t_0  \rfloor +1$ times, we obtain the following near-diagonal lower bound
\begin{equation}\label{e:main4.2}
p^{\kappa}(t,x,y)\ge C_3 \pinv(t)^{-d} \ \ \textrm{for all }t\in(0,T] \textrm{ and }|x-y|\le C_4\pinv(t)\, 
\end{equation}
for some constants $C_3,C_4>0$. Indeed, for $t \in (0,T]$ and $x,y \in \R^d$ with $|x-y| \le C_4 \pinv(t)$ where $C_4:=C_2(n_0/a_1)^{1/\alpha_1} \in (0, 2^{-4})$ is a sufficiently small constant satisfying
$$ C_4 \pinv(t) \le C_2 \pinv(\frac{t}{n_0}) $$
by \eqref{e:wsc-inv}. Let $n = \lfloor \frac{T|x-y|}{C_4t_0  \pinv(t)} \rfloor +1 $ and $z_1, \dots, z_{n-1}$ be the points in the segment from $x$ to $y$ satisfying $|z_1-x| = |z_{i+1} - z_i | = |y - z_{n-1}| = \frac{|x-y|}{n}$. Note that $n \le n_0$ since $|x-y| \le C_4 \pinv(t)$. Then, by \eqref{e:main4-1} 
\begin{align*}
&p^k(t,x,y) = \int_{\R^d} \dots \int_{\R^d} p^\kappa(t/n,x,w_1) p^\kappa(t/n,w_1,w_2) \dots p^\kappa(t/n,w_{n-1},y)        dw_1 \dots dw_{n-1} \\
&\ge \int_{B(z_1,\frac{|x-y|}{3n})} \dots \int_{B(z_{n-1}, \frac{|x-y|}{3n})}p^\kappa(t/n,x,w_1) p^\kappa(t/n,w_1,w_2) \dots p^\kappa(t/n,w_{n-1},y)        dw_1 \dots dw_{n-1} \\
&\ge \omega_d^{n-1} \big(\frac{|x-y|}{3n} \big)^{-(n-1)d} \cdot \frac{C_1^n}{2^n} \pinv(\frac{t}{n})^{-nd} \ge C_3 \pinv(t)^{-d}.
\end{align*}
Here we used $\frac{|x-y|}{n} \ge  \frac{T}{C_4 t_0} \pinv(t)$, $n \le n_0$ and \eqref{e:wsc-inv} for the last line. Therefore, we obtain \eqref{e:main4.2}. 
Now we assume $|x-y| > C_4\pinv(t)$ and let $ \lam>0  $ be the constant in Lemma \ref{l:ep} for $\eps=1/2$
and $\tau(z,c, t):=\tau_{B(z,c\pinv(t))}$. Then for every $0<t \le T$,
\begin{equation}\label{e:main4.3}
\sup_{z \in \R^d} \P_z\big(\tau(z,2^{-2}C_4, t)
\le \lambda t\big)\le \frac12\, .
\end{equation}
Let $\sigma=\inf\{s\ge 0: X_s\in B(y,2^{-2}C_4\pinv(t))\}$ be the hitting time of $B(y,2^{-2}C_4\pinv(t))$.
By the strong Markov property and \eqref{e:main4.3} we have
\begin{align}
\begin{split}\label{e:main4.4}
&\P_x\Bigg(X_{\lambda t}\in B(y,2^{-1}C_4\pinv(t))\Bigg) \ge \P_x\left(\sigma\le \lambda t, \sup_{s\in [\sigma, \lambda t]} |X_s-X_{\sigma}|<2^{-2}C_4\pinv(t)\right)\\
&=\E_x\left(\P_{X_{\sigma}}\Big(\sup_{s\in[0,\lambda t]} |X_s-X_0|<2^{-2}C_4\pinv(t)\Big); \sigma\le \lambda t\right)\\
&\ge  \inf_{z\in  B(y,\pinv(t))} \P_z\big(\tau(z,2^{-2}C_4, t)>\lambda t\big) \P_x\big(\sigma \le \lambda t\big)\\
&\ge \frac12 \P_x\big(\sigma \le \lambda t\big) \, \ge \, \frac12 \P_x\left(X_{\lambda t \wedge \tau(x,2^{-3}C_4, t)}\in B(y,2^{-2}C_4\pinv(t))\right)\, .
\end{split}
\end{align}
Since $|x-y| > C_4\pinv(t)$, we have 
$$
X_s\notin B\left(y,2^{-2}C_4\pinv(t)\right)^c \subset B\left(x,2^{-3}C_4\pinv(t)\right)^c\,, \quad s<\lambda t \wedge \tau(x,2^{-3}C_4, t).
$$
Thus,
$$
{\mathbf 1}_{\{X_{\lambda t \wedge \tau(x,2^{-3}C_4, t)}\in B(y,2^{-2}C_4\pinv(t))\}}=\sum_{s\le \lambda t \wedge \tau(x,2^{-3}C_4, t)} {\mathbf 1}_{\{X_s\in B(y,2^{-2}C_4\pinv(t))\}}\, .
$$
Therefore, by the L\'evy system formula in \eqref{e:LSF} we obtain
\begin{align}
&\P_x\left(X_{\lambda t \wedge \tau(x,2^{-3}C_4,t)}\in B(y,2^{-2}C_4\pinv(t))\right) 
=\ \E_x\left[\int_0^{\lambda t \wedge \tau(x,2^{-3}C_4,t) }\int_{B(y,2^{-2}C_4\pinv(t))}J_X(X_s,u)\, du\, ds\right] \nn \\
&\ge\ \E_x\left[\int_0^{\lambda t \wedge \tau(x,2^{-3}C_4,t) }\int_{B(y,2^{-2}C_4\pinv(t))}\kappa_0 J(|X_s-u|)
{\bf 1}_{\{u:|X_s-u| \le |x-y|\}}
\, du\, ds\right]. \label{e:main4.5}
\end{align}
Let $w$ be the point in the segment from $x$ to $y$ satisfying $|w-y| = 3 \cdot 2^{-4} C_4 \pinv(t) $. Since $|x-X_s| \le 2^{-3}C_4 \pinv(t)$, we have that for any $u \in B(w, 2^{-4}C_4\pinv(t) )$,
\begin{equation*}
|X_s - u| \le |x-X_s| + |x-w| + |w-u| \le |x-y|
\end{equation*}
so that $B(w, 2^{-4}C_4 \pinv(t)) \subset B(y,2^{-2}C_4 \pinv(t)) \cap \{u:|X_s-u| \le |x-y|\}$ for every $s<\lambda t \wedge \tau(x,2^{-3}C_4, t).$ Thus, 
\begin{align}
\begin{split}\label{e:main4.6}
&\E_x\left[\int_0^{\lambda t \wedge \tau(x,2^{-3}C_4,t)}\int_{B(y,\pinv(t))}\kappa_0 J(|X_s-u|)
{\bf 1}_{\{u:|X_s-u| \le |x-y|\}}
\, du\, ds\right] \\
&\ge \kappa_0 \E_x\left[\int_0^{\lambda t \wedge \tau(x,2^{-3}C_4,t)}\int_{B(w,2^{-4}C_4 \pinv(t))} J(|x-y|)
\, du\, ds\right] \\
&\ge c_1 \pinv(t)^{d} J(|x-y|) \E_x [ \lam t \land \tau(x,2^{-3}C_4, t) ] \\
&\ge c_2 t \pinv(t)^{d} J(|x-y|),
\end{split}
\end{align}
where we used \eqref{e:exit} and \eqref{e:wsc-Phi} for the last line.

Therefore, combining \eqref{e:main4.4}, \eqref{e:main4.5} and \eqref{e:main4.6} we arrive
\begin{align*}
	p^{\kappa}(t,x,y)&\ge  \int_{B(y,\pinv(t))}p^{\kappa}(\lambda t, x,z)p^{\kappa}((1-\lambda)t,z,y)\, dz\\
	&\ge \inf_{z\in B(y,2^{-1}C_4\pinv(t))}p^{\kappa}((1-\lambda)t,z,y) \int_{B(y,2^{-1}C_4\pinv(t))}p^{\kappa}(\lambda t, x,z)\, dz\\
	&\ge c_3  t\pinv(t)^{-d} \pinv(t)^d J(|x-y|).
\end{align*}
for all $0<t \le T$ and $x,y \in \R^d$ with $|x-y| > C_4\pinv(t)$.
\qed

\textbf{Proof of Theorem \ref{t:main1}.} 
By Remarks \ref{r:KSV} and \ref{r:KSV1}, 
$p^\kappa(t,x,y)$ defined in \eqref{d:p^kap} satisfies \eqref{e:main1-1}, \eqref{e:main1-3} and \eqref{e:main1-4}. Also, \eqref{e:main1-2} and \eqref{e:main4-1} follow from Theorem  \ref{t:p^kap} and Lemma \ref{l:main2}, respectively. 
It remains to  show the uniqueness part of Theorem \ref{t:main1}.
Recall that we observe in Remark \ref{r:KSV1} that \cite[(1.9)]{KSV} holds. Thus all results in  \cite[Sections 5.1 and 5.2]{KSV} hold for our case.
Since properties (i)--(iii) are stronger than ones in  \cite[Theorem 1.1]{KSV}, we now see that 
the proof of  the uniqueness part of Theorem \ref{t:main1} is exactly same as the one of the uniqueness part of \cite[Theorem 1.1]{KSV}. 
\qed

\bigskip
\noindent
{\bf Acknowledgements:} After finishing the main part of this paper, the result was announced by the second named author  during 
{\it Extended Kansai Probability Seminar} on January 27, 2017. We thanks Research Institute for Mathematical Sciences and the Department of Mathematics
of Kyoto University for the hospitality.
\end{doublespace}

\bibliography{Exp}{}
\bibliographystyle{alpha}

\vskip 0.1truein

\parindent=0em

{\bf Panki Kim}

Department of Mathematical Sciences and Research Institute of Mathematics,

Seoul National University, Building 27, 1 Gwanak-ro, Gwanak-gu Seoul 08826, Republic of Korea

E-mail: \texttt{pkim@snu.ac.kr}

\bigskip

{\bf Jaehun Lee}

Department of Mathematical Sciences,

Seoul National University, Building 27, 1 Gwanak-ro, Gwanak-gu Seoul 08826, Republic of Korea

E-mail: \texttt{hun618@snu.ac.kr}

\end{document}